\documentclass[preprint]{elsarticle}
\usepackage[utf8]{inputenc}  
\usepackage[T1]{fontenc} 
\usepackage{lmodern}
\usepackage{amsmath} 
\usepackage{amssymb}
\usepackage{amsthm} 
\usepackage{mathrsfs}
\usepackage{graphicx}
\usepackage{ccaption}
\usepackage{hyperref}
\usepackage{listings}
\usepackage{enumitem}
\usepackage{tikz}
\usepackage{bm}
\usepackage{float}
\usepackage{pgfplots}
\pgfplotsset{compat=newest}
\usetikzlibrary{shapes}
\usetikzlibrary{patterns}
\usepackage{stmaryrd}
\usepackage[pagewise]{lineno}

\newdefinition{defi}{Definition}
\newdefinition{asum}{Assumption}
\newtheorem{theorem}{Theorem}
\newtheorem{prop}{Proposition}
\newtheorem{coro}{Corollary}
\newtheorem{lem}{Lemma}
\newdefinition{rem}{Remark}
\newcommand{\R}{\mathbb{R}}
\newcommand{\N}{\mathbb{N}}

\newcommand{\C}{\mathbb{C}}
\newcommand{\dd}{\mathrm{d}}

\newcommand{\eps}{\varepsilon}
\newcommand{\loc}{\text{loc}}

\newcommand{\Ai}{\mathcal{A}}
\newcommand{\Bi}{\mathcal{B}}

    \pgfplotsset{
        colormap={parula}{
            rgb255=(53,42,135)
            rgb255=(15,92,221)
            rgb255=(18,125,216)
            rgb255=(7,156,207)
            rgb255=(21,177,180)
            rgb255=(89,189,140)
            rgb255=(165,190,107)
            rgb255=(225,185,82)
            rgb255=(252,206,46)
            rgb255=(249,251,14)
        },
    }
\everymath{\displaystyle}   
\newcommand\tsum{\textstyle\sum\nolimits}
\newcommand\tint{\textstyle\int\nolimits}

\def\wt{\widetilde}
\def\ph{\varphi}
\def\top{\text{top}}
\def\bot{\text{bot}}
\def\cH{\mathcal{H}}
\def\cO{\mathcal{O}}
\def\tH{\text{H}}
\def\tL{\text{L}}
\DeclareMathOperator\supp{\text{supp}}

\begin{document}

\begin{frontmatter}
\title{The Helmholtz problem in slowly varying waveguides at locally resonant frequencies}
\author[1]{\'Eric Bonnetier}
\ead{Eric.Bonnetier@univ-grenoble-alpes.fr}
\author[2]{Angèle Niclas \corref{cor1}}
\ead{angele.niclas@ec-lyon.fr}
\author[2]{Laurent Seppecher}
\ead{laurent.seppecher@ec-lyon.fr}
\author[2]{Grégory Vial}
\ead{gregory.vial@ec-lyon.fr}
\cortext[cor1]{Corresponding author}
\affiliation[1]{organization={Institut Fourier, Université Grenoble Alpes}, adressline={100 Rue des Mathématiques}, postcode={38610}, city={Gières}, country={France}}
\affiliation[2]{organization={Institut Camille Jordan, \'Ecole Centrale Lyon}, adressline={36 Avenue Guy de Collongue}, postcode={69134}, city={Écully}, country={France}}
\begin{abstract}
This article aims to present a general study of the Helmholtz problem in slowly varying waveguides. This work is of particular interest at locally resonant frequencies, where a phenomenon close to the tunnel effect for Schrödinger equation in quantum mechanics can be observed.  In this situation, locally resonant modes propagate in the waveguide under the form of Airy functions. Using previous mathematical results on the Schrödinger equation, we prove the existence of a unique solution to the Helmholtz source problem with outgoing conditions in such waveguides. We provide an explicit modal approximation of this solution, as well as a control of the approximation error in $\text{H}^1_\loc$. The main theorem is proved in the case of a waveguide with a monotonously varying profile and then generalized using a matching strategy. We finally validate the modal approximation by comparing it to numerical solutions based on the finite element method. 
\end{abstract}

\end{frontmatter}

\section{Introduction}

In this article, we study the propagation in the time harmonic regime of waves  generated by sources in a slowly varying waveguide of dimension 2. The waveguide is described by
\begin{equation}\label{eq:omegatilde}
\widetilde{\Omega}:=\left\{(x,y)\in \R^2\ |\ 0< y < h(x)\right\},
\end{equation}
where $h\in \mathcal{C}^2(\R)\cap W^{2,\infty}(\R)$ is a positive profile function defining the top boundary. Here, the bottom boundary is assumed to be flat but a similar analysis can be conducted with both slowly varying top and bottom boundaries. In the time harmonic regime, the wave field $\widetilde{u}$ satisfies the Helmholtz equation with Neumann boundary conditions 
\begin{equation}\label{eqdebut} \quad 
\left\{\begin{array}{cl} \Delta \widetilde{u}+k^2\widetilde{u} =-\wt f & \text{ in } \widetilde{\Omega},\\
 \partial_\nu \widetilde{u} = \wt b & \text{ on } \partial\widetilde{\Omega}, \end{array}\right. 
\end{equation} 
where $k>0$ is the frequency, $\wt f$ is an interior source term and $\wt b$ is a possible boundary source term. In this work, a waveguide is said to be slowly varying when there exists a small parameter $\eta>0$ such that $\Vert h'\Vert_{\text{L}^\infty(\R)}\leq \eta$ and $\Vert h''\Vert_{\text{L}^\infty(\R)}\leq \eta^2$. Such waveguides are good models of ducts or corroded pipes, and studying the sound transmission through this type of structure can be used to reduce noise emission (see \cite{nielsen1}) or to perform non destructive monitoring of pipes or blood vessels (see \cite{honarvar1}).

\subsection{Scientific context}

Wave propagation in varying waveguides, whether acoustic or elastic, has already been studied by several authors. From a numerical point of view, the articles \cite{lu2, fabro1, mitsoudis1} give different methods to adapt the finite element method to numerically compute the wave field in varying waveguides. In \cite{pagneux2}, the authors study from the theoretical point of view the propagation of waves in a general varying elastic waveguide using a modal decomposition. The same kind of method is used in \cite{folguera1}, in the case of a slowly varying waveguide. However, in these articles, the authors choose to avoid all the locally resonant frequencies of the waveguide, which are the frequencies $k>0$ such that $k=\pi n/h(x^\star)$ for a mode $n\in \N$ and a longitudinal position $x^\star\in \R$.

In another approach, the authors of \cite{perel1,nielsen1,galanenko1} choose to work near locally resonant frequencies of the waveguide. They mainly show that this problem is very close to the tunneling effect seen in quantum mechanics for the Schrödinger equation (see for instance \cite{roy1}). Indeed, the wave field can be decomposed as a sum of modes $\widetilde{u}(x,y)=\tsum_{n\in \N} u_n(x) \varphi_n(y)$ (see section 2 for more details) and when $k=n\pi/h(x^\star)$, for some $x^\star\in\R$, the equation satisfied by the mode $u_n$ is close to the Schrödinger equation 
\begin{equation}
\partial_{xx} u_n(x)+(V(x)-E)u_n(x)=0,
\end{equation}
where $V$ and $E$ depend on $h$, $\wt f$ and $\wt b$ and satisfy $V(x^\star)-E=0$. This equation for a simple mode was studied from a mathematical point of view by F. W. J. Olver in \cite{olver1} and \cite{olver2}, and it was proved that the solution $u_n$ could be expressed using Airy functions of the first and second kind \cite{abramowitz1}. In all the articles \cite{perel1,nielsen1,galanenko1}, the same methodology is used: firstly, the authors assume that mode coupling is negligible in a slowly varying waveguide. Under this so-called adiabatic approximation, every mode is independent from the others. Secondly, they seek solutions expressed as Wentzel–Kramers–Brillouin (WKB) asymptotic series (see \cite{olver3}), and they use the study of the Schrödinger equation to find an approximation of the wave field in the waveguide. 

Our work is inspired by this methodology and provides a similar approximation of the wave field in a slowly varying waveguide. However, contrary to the work mentioned above, we are not making any \textit{a priori} assumptions on the wave field such as WKB asymptotic development or the adiabatic decoupling of modes. Like them, we use the work of \cite{olver1,olver2} to get an approximation of the Schrödinger equation, but we improve it by providing precise control of the approximation error. More importantly, we provide a way to justify the adiabatic decoupling of the modes using a Born approximation of the wave field. To this end, we again rely on \cite{olver1,olver2} to control the wave field by the general source term that generated it. The main result of our article is given by Theorem \ref{th1} that proves the existence of a unique solution of the problem \eqref{eqdebut} when $h$ is an increasing function in $\mathcal{C}^2(\R)\cap W^{2,\infty}(\R)$ and when $\eta$ is small enough (compared to $\text{supp}(h')$, $\min(h)$ and the distance between $k$ and the left and right resonances of the waveguide). This theorem also provides an approximation of the solution of \eqref{eqdebut} and a control of the approximation error. 

Finally, we provide a numerical validation of the approximation of the wave field in a slowly varying waveguide. By comparing our approximation to solutions generated by a finite element method, we show that this approximation is an excellent tool to  numerically compute the wave field in a slowly varying waveguide in a very fast way.

\subsection{Outline of the paper} 
The paper is organized as follows. In section 2, we briefly explain the modal decomposition in general waveguides, and we recall classical results used in the rest of the paper. In section 3, we study the particular case of a slowly varying waveguide where the width $h$ is an increasing function of $x$, and we prove Theorem 1. In section 4, we adapt the method developed in section 3 to the general case of a varying waveguide provided the variations of the profile are sufficiently slow, and we describe more precisely the cases of compressed or dilated waveguides. In section 5, we numerically validate our results by comparing the approximations derived in sections 3 and 4 with the solutions generated using a finite element solver with PML (perfectly matched layers, see \cite{berenger1}) in a truncated waveguide.

\subsection{Notations}

The varying waveguide is denoted by $\widetilde{\Omega}$, and its boundary $\partial\widetilde{\Omega}$. The subscript “top” (resp. “bot”) indicates the upper boundary of the waveguide (resp. lower). The straight waveguide is defined by $\Omega=\R\times (0,1)$, and its boundary is denoted by $\partial\Omega$. For every $r>0$, we set $\Omega_r=(-r,r)\times (0,1)$ and $\Gamma_r=(-r,r)\times\{0\} \cup (-r,r)\times \{1\}$. For both, $\nu$ denotes the outer normal vector. The spaces $\tH^1$, $\tH^2$, $W^{1,1}$ $\tH^{1/2}$ over $\wt \Omega$, $\Omega$ or their boundaries are the usual Sobolev spaces on piece wise smooth domains. The space $\widetilde{\text{H}}^{1/2}(-r,r)$ is the closure of $\mathcal{D}(-r,r)$, the space of distributions with support in $(-r,r)$, for the $\text{H}^{1/2}(\R)$ norm (see \cite{mclean1} for more details).

The operator norm between two Banach spaces $E_1$ and $E_2$ is denoted $\Vert \cdot \Vert_{E_1, E_2}$ and is defined for every linear operator $\mathcal{S} : E_1 \rightarrow E_2$ by 
\begin{equation}
\Vert \mathcal{S} \Vert_{E_1,E_2}:=\sup_{x\in E_1 \, |\,\Vert x\Vert_{E_1}=1} \Vert \mathcal{S}(x) \Vert_{E_2}.
\end{equation}

The Airy function of the first kind (resp. second kind) is denoted by $\Ai$ (resp. $\Bi$). These functions are linear independent solutions of the Airy equation $y''-xy=0$ (see \cite{abramowitz1} for more results about Airy functions) and are depicted in Figure \ref{airy}. 

\begin{figure}[h!]
\begin{center}
\scalebox{.5}{\input{airy}}
\caption{\label{airy} Representation of the Airy functions $\Ai$ and $\Bi$.}
\end{center}
\end{figure}

\section{Modal decomposition and local wavenumbers in a varying waveguide}

In this section, we recall some classical results about modal decompositions, the proofs of which can be found in \cite{bonnetier1, bourgeois1}.

\begin{defi}\label{def:modes}  We define the sequence of functions $(\wt \ph_n)_{n\in\N}$
\begin{equation}\label{phintilde}
\forall (x,y)\in\wt \Omega,\quad\wt \ph_n(x,y) :=
\left\{\begin{array}{cl}
1/\sqrt{h(x)}\quad  &\text{if } n=0, \\ 
 \frac{\sqrt{2}}{\sqrt{h(x)}}\cos\left(\frac{n\pi y}{h(x)}\right)\quad &\text{if } n\geq 1,
\end{array}\right.
\end{equation}
which for any $x\in \R$ defines an orthonormal basis of $\tL^2(0,h(x))$. In the special case of a regular waveguide where $h=1$ everywhere, this sequence of functions is independent of $x$, takes the form 
\begin{equation}\label{phin}
\forall y\in(0,1),\quad \varphi_n(y) :=
\left\{\begin{array}{cl}
1\quad  &\text{if } n=0,\\ 
\sqrt{2}\cos\left(n\pi y\right)\quad &\text{if } n\geq 1.
\end{array}\right.
\end{equation}
and defines an orthonormal basis of $\tL^2(0,1)$.
\end{defi}

Hence, any solution $\wt u \in \tH^2_{\loc}\big(\wt \Omega\big)$ of \eqref{eqdebut} admits a unique modal decomposition
\begin{equation}\label{decmode}
\wt u(x,y)=\sum_{n\in \N}\wt u_n(x)\wt \ph_n(x,y)\quad\text{where}\quad \wt u_n(x):=\int_0^{h(x)}\wt u(x,y)\wt \ph_n(x,y)\dd y.
\end{equation}
Note that when $h$ is constant (outside of $\supp h'$), each mode $\wt u_n$ satisfies the simple equation $\wt u_n''+k_n^2\wt u_n=-\wt g_n$ where $k_n$ is the wavenumber. When $h$ is variable, the decomposition \eqref{decmode} motivates the following definition:

\begin{defi}\label{def:wavenumber}  The local wavenumber function of the mode $n\in\N$ is the complex function $k_n:\R\to \C$ defined by
\begin{equation}\label{kn}
k_n^2(x):=k^2-\frac{n^2\pi^2}{h(x)^2},
\end{equation}
with $\text{Re}(k_n), \text{Im}(k_n)\geq 0$. 
\end{defi}

One of the main difficulties of this work is that as $h(x)$ is non constant, $k_n(x)$ can vanish for some $x\in\R$ and change from a positive real number to a purely imaginary number. We distinguish three different situations. 

\begin{defi}  A mode $n\in\N$ falls in one of the three following situations:
\begin{enumerate}
\item If $n>kh(x)/\pi$ for all $x\in\R$ then $k_n(x)\in(0,+\infty)$ for all $x\in\R$ and the mode $n$ is called propagative. 
\item If $n<kh(x)/\pi$ for all $x\in\R$ then $k_n(x)\in i(0,+\infty)$ for all $x\in\R$ and the mode $n$ is called evanescent. 
\item If there exists $x^\star\in \R$ such that $n=kh(x^\star)/\pi$ the mode $n$ is called locally resonant. The associated points $x^\star$ are called resonant points. They are simple if $h'(x^\star)\neq 0$, and multiple otherwise. 
\end{enumerate}
A frequency $k>0$ for which there exists at least a locally resonant mode is called a locally resonant frequency.  
\end{defi}

Using the wavenumber function, one can adapt the classic Sommerfeld (or outgoing) condition, defined in \cite{bonnetier1} for regular waveguides, to general varying waveguides $\wt\Omega$. This condition will be used later to guarantee uniqueness for the source problem \eqref{eqdebut}.

\begin{defi}\label{def:outgoing}  A wavefield $\wt u \in \tH^2_{\loc}\big(\wt \Omega\big)$ is said to be outgoing if it satisfies 
\begin{equation} \label{sommer}\left| \wt u_n'(x)\frac{x}{|x|}-ik_n(x)\wt u_n(x) \right| \underset{|x|\rightarrow +\infty}{\longrightarrow} 0 \qquad \forall n\in \N,
\end{equation}
where $\wt u_n$ is given in \eqref{decmode}. 
\end{defi}

\section{The Helmholtz equation in a waveguide with increasing width}

In all this work, we make the following assumptions:

\begin{asum}\label{def:slow} We assume that $h\in \mathcal{C}^2(\R)\cap W^{2,\infty}(\R)$ and satisfies 
\begin{equation}\nonumber
\forall x\in \R \quad 0<h_{\min} \leq h(x) \leq h_{\max} <\infty,
\end{equation}
\begin{equation}\nonumber
\Vert h'\Vert_{\text{L}^\infty(\R)} <\eta, \quad \Vert h''\Vert_{\text{L}^\infty(\R)}<\eta^2, \quad \supp h'\subset \left(-\frac{R}{\eta},\frac{R}{\eta}\right),
\end{equation}
for some $\eta>0$ and $R>0$. 
\end{asum}

\begin{figure}[h]
\begin{center}
\begin{tikzpicture} 
\draw (-1.5,0) -- (9,0);
\draw (-1.5,1) -- (0,1); 
\draw (7,1.393) -- (9,1.393); 
\draw (4.5,0.6) node{$\widetilde{\Omega}$}; 
\draw [domain=0:7, samples=100] plot (\x,{1+12/7/7/7/7/7*\x*\x*\x*(\x*\x/5-7*\x/2+49/3)});  
\end{tikzpicture}
\end{center}
\caption{\label{fig1} A waveguide with increasing width.}
\end{figure}
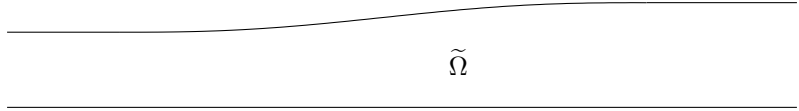

Moreover, in all this section, and we assume that $h$ is increasing in $\supp h'$ (the general case will be treated in section \ref{section4}). Such a waveguide is represented in Figure \ref{fig1}. 

The aim of this section is to state and prove the main theorem of this work, which is a triple result. It provides existence and uniqueness of the solution $\wt u\in \tH^2_{\loc}\big(\wt \Omega\big)$ to the source problem
\begin{equation} \label{MAFT} \quad  
\left\{\begin{array}{cl} \Delta \widetilde{u} +k^2 \widetilde{u} =-\wt f  &\text{ in } \widetilde{\Omega}, \\
\partial_\nu \widetilde{u} = \wt b_\top  &\text{ on } \partial\widetilde{\Omega}_\text{top}, \\
\partial_\nu \widetilde{u} = \wt b_\bot  &\text{ on } \partial\widetilde{\Omega}_\text{bot},\\
\widetilde{u} \text{ is outgoing}. \end{array}\right. \tag{$\widetilde{\mathcal{H}}$}
\end{equation}
It also gives an explicit modal approximation of $\wt u$, and provides a computable error bound for the local $\tH^1$-norm between the approximation and the exact solution. We first explain our strategy for proving such an existence result.


\subsection{Sketch of proof}

In order to use results on the modal decomposition in the regular waveguide, the first step is to map the perturbed waveguide to the regular one using the canonical mapping $\psi:\Omega\to\wt \Omega$ defined by $\psi : (x,y) \mapsto (x, h(x)y)$. The problem $\big(\wt\cH\big)$ is then equivalent to
\begin{equation} \label{eqdif}
\left\{\begin{array}{cl} \Delta_hu +k^2 u =-f & \text{ in } \Omega, \\
\partial_\nu u -D_h u =  b_\top & \text{ on } \partial\Omega_\text{top}, \\
\partial_\nu u = b_\bot & \text{ on } \partial\Omega_\text{bot},\\ u \text{ is outgoing}, \end{array}\right. \tag{$\mathcal{H}$}
\end{equation}
where $u:=\wt u\circ\psi$, $f:=\wt f\circ\psi$, $b_\top:=\wt b_\top \sqrt{1+(h')^2}/h$ and $b_\bot = \wt b_\bot/h$. The operators $\Delta_h$ and $D_h$ are differential operators of order two and one respectively (see their expressions in \eqref{eqomega}). As there is no easy way to solve explicitly this equation we shall approach it by a simpler problem. To this end, we neglect the small terms in the operators $\Delta_h$ and $D_h$, which depend on $h'$ and $h''$. This leads to a much simpler problem that reads 
\begin{equation}\label{eq1}
\left\{\begin{array}{cl}\displaystyle \partial_{xx} v +\frac{1}{h(x)^2}\partial_{yy} v+k^2 v = -f & \text{ in } \Omega, \\
\partial_\nu v = b_\top & \text{ on } \partial\Omega_\text{top},\\
\partial_\nu v = b_\bot & \text{ on } \partial\Omega_\text{bot}, \\
v \text{ is outgoing}.
 \end{array}\right. \tag{$\mathcal{H}'$}
\end{equation}

Next, we seek a solution to \eqref{eq1} in a modal form. To this end, we use the fact that the wave field $v$ and the source $f$ can be decomposed in a sum of modes in the straight guide

\begin{equation}\label{eq:decomp}
\begin{aligned}
v(x,y)=\sum_{n\in \N}v_n(x)\ph_n(y)\quad\text{where}\quad v_n(x):=\int_0^1w(x,y) \ph_n(y)\dd y,\\
f(x,y)=\sum_{n\in \N}f_n(x)\ph_n(y)\quad\text{where}\quad f_n(x):=\int_0^1f(x,y) \ph_n(y)\dd y.
\end{aligned}
\end{equation}
We deduce that the two dimensional problem $(\cH')$ is equivalent to a sequence of one dimensional problems for $n\in \N$: 
\begin{equation}\label{multi1D}
 \left\{\begin{array}{cl} v_n''(x) +k_n(x)^2 v_n(x) =-(f_n+\ph_n(1)b_\top+ \ph_n(0)b_\bot)(x) & \text{ in } \R, \\ v_n \text{ is outgoing}, & \end{array}\right. \tag{$\mathcal{H}'_n$}
\end{equation}
where $k_n$ is the local wavenumber function given in \eqref{kn}. As explained in \cite{bourgeois1}, the modal problem $(\cH_n')$ is well-defined when $\{x\in \R \,|\, k_n(x)=0\}$ has zero measure. Since $h$ is increasing on $\supp h'$, this only occurs when $k\neq n\pi/h_{\min}$ or $k\neq n\pi/h_{\max}$. We assume this is the case and that 
\begin{equation}
\delta=\min_{n\in \N}\left(\sqrt{\left|k^2-\frac{n^2\pi^2}{{h_{\min}}^2}\right|},\sqrt{\left|k^2-\frac{n^2\pi^2}{{h_{\max}}^2}\right|} \right)>0.
\end{equation}
which is supposed to be positive.

For each type of modes $n\in \N$, the study of the equation $v_n''+k_n^2 v_n=0$ has been carried out in \cite{olver1} and \cite{olver2}. We summarize the main ideas bellow. The analysis depends on whether the mode $n$ is propagative, evanescent or locally resonant.

\begin{itemize}
\item[(a)] If $n$ is propagative or evanescent, $|k_n(x)|>0$ for all $x\in \R$ and we set $z(x)=\tint^x |k_n|$. There is has a one-to-one correspondence between $x$ and $z$ and if we define $w_n=\sqrt{z}v_n$, then $w_n$ satisfies the partial differential equation
\begin{equation}
\partial_{zz}w_n  \pm w_n= \zeta(x,z)\,w_n,
\end{equation}
where $\Vert \zeta\Vert_{L^\infty(\R)} =\cO(\eta)$. The solutions of $y''=\pm y$ are exponential functions, and since $\zeta$ is small, we can prove that $w_n$ is almost equal to a sum of two exponential functions and we can control the approximation error (see  \eqref{olv1} and \eqref{olv2} for more details).

\item[(b)] If $n$ is locally resonant, since $h$ is increasing on $\supp h$, there is a single resonant point $x^\star\in\R$ and we define 
\begin{equation}\label{eq:xi}
\xi(x):=
\left\{\begin{aligned}
\left(-\frac{3}{2}i\int_x^{x^\star}k_n(t)\dd t\right)^{2/3} & \text{ if } x<x^\star, \\ -\left(\frac{3}{2}\int_{x^\star}^x k_n(t) \dd t \right)^{2/3} & \text{ if } x>x^\star.
\end{aligned}\right.\end{equation}
This new variable is in one-to-one correspondence with $x$ and if we denote $w_n=-(\sqrt{k_n}/\xi^{1/4})v_n$, then $w_n$ satisfies the partial differential equation 
\begin{equation}
\partial_{\xi\xi} w_n - \xi w_n= \zeta(\xi)\,w_n,
\end{equation}
where $\Vert \zeta\Vert_{L^\infty(\R)} =\cO( \eta)$. The solutions of $y''= x y$ are known as the Airy functions, and since $\zeta$ is small, we can prove that $w_n$ is approximated by a sum of Airy functions and we can control the approximation error (see  \eqref{olv3} for more details). 

\end{itemize}

Using these results, we prove that $(\cH_n')$ has a unique solution and we provide an explicit approximation of this solution. By equivalence, this approach yields the unique solution to $(\cH')$ and its approximation. With a control of the approximation error between $(\cH)$ and $(\cH')$, we obtain an explicit approximation of $(\cH)$ and by change of variable, of $(\mathcal{\widetilde{H}})$.

\subsection{Main result}

We now state the main result of this work, which shows existence and uniqueness of the solution $u$ of \eqref{eqdif} (and thus of the solution of \eqref{MAFT}) and provides an approximation of $u$ with control of the approximation error in $\text{H}^1_\loc$.

\begin{theorem}\label{th1}
Let $h$ be an increasing function which defines a varying waveguide $\tilde \Omega$ that satisfies assumption \ref{def:slow}. Consider sources $f\in \tL^2(\Omega)$, $\bm{b}=(b_\bot,b_\top)\in ({\text{H}}^{1/2}(\R))^2$ both with compact support contained in $\Omega_r$ and $\Gamma_r$ respectively, for some $r>0$. Assume that there is a unique locally resonant mode $N\in \N$, associated with a simple resonant point $x^\star\in\R$. 

There exists $\eta_0>0$, depending only on $h_{\min}$, $h_{\max}$, $\delta$, $r$ and $R$, such that if $\eta<\eta_0$, then the problem $(\cH)$ admits a unique solution $u\in \text{H}^2_{\text{loc}}\big(\Omega)$. Moreover, this solution is approximated by $u^{\text{app}}$ defined for almost every $(x,y)\in \Omega$ by 
\begin{equation}\label{greentot}
u^{\text{app}}(x,y)=\sum_{n\in \N} \int_\R G_n^{\text{app}}(x,s)\big(f_n+\ph_n(1)b_\top+\ph_n(0)b_\bot\big)(s)\dd s \, \ph_n\left(y\right),
\end{equation}
where $f_n$ is defined in \eqref{eq:decomp}, $\ph_n$ is defined in \eqref{phin} and $G_n^{\text{app}}(x,s)$ is equal to 
\begin{equation}\label{greenfunction}
\left\{\begin{aligned}
&\frac{i}{2\sqrt{k_n(s)k_n(x)}}\exp\left(i\left|\int_s^xk_n\right|\right), & \,\,\text{ if } n<N,\\
&\frac{1}{2\sqrt{|k_n|(s)|k_n|(x)}}\exp\left(-\left|\int_s^x|k_n|\right|\right), &\,\,\text{ if } n>N,\\
&\left\{\begin{aligned}
\frac{\pi(\xi(s)\xi(x))^{1/4}}{\sqrt{k_n(s)k_n(x)}}\big(i\Ai+\Bi\big)\circ\xi(s)\Ai\circ\xi(x)& \,\,\text{ if } x<s, \\
\frac{\pi(\xi(s)\xi(x))^{1/4}}{\sqrt{k_n(s)k_n(x)}}\big(i\Ai+\Bi\big)\circ\xi(x)\Ai\circ\xi(s)& \,\,\text{ if } x>s, \\
\end{aligned}\right.
&\,\,\text{ if } n=N.
\end{aligned}\right.\end{equation}
Th function $k_n$ is the wavenumber function defined in definition \ref{def:wavenumber} and the function $\xi$ is given in equation \eqref{eq:xi}. Moreover, there exists a constant $C>0$ depending only on $h_{\min}$, $h_{\max}$, $\delta$, $r$ and $R$ such that
\begin{equation}
\Vert u -u^{\text{app}}\Vert_{\text{H}^1(\Omega_r)}\leq \eta C \left(\Vert f\Vert_{\text{L}^2(\Omega)}+\Vert \bm{b}\Vert_{\left({\text{H}}^{1/2}(\R)\right)^2}\right).
\end{equation}

\end{theorem}

\begin{rem}
If there are no resonant modes, the result can be adapted by deleting the line $n=N$ in \eqref{greenfunction}. On the other hand, if there are multiple locally resonant modes, the third line of \eqref{greenfunction} becomes true for every resonant mode. 
\end{rem}

\begin{rem}
If $\wt\Omega$ is a regular waveguide, we find the same expression for the wave field as in \cite{bonnetier1}. We also see that the behavior of propagative and evanescent modes in a perturbed waveguide is similar to that in a regular waveguide. The term $\tint_s^x|k_n|$ simply acts as a change of variable in the phase. 
\end{rem}

\begin{rem}
Looking at the proof, we can see that the constant $C$ has a dependence on $\delta$, $r$ and $R$ of the form $C=\cO(r^2\delta^{-6}+rR\delta^{-8})$. Doing the same proof using $\text{W}^{2,\infty}$ spaces instead of $\text{H}^2$, we can also prove that for every $x\in \R$,  
\begin{equation}
|u_{N}(x)-u_N^{\text{app}}(x)|\leq \eta C \left(\Vert f\Vert_{\text{L}^\infty(\R)}+\Vert \bm{b}\Vert_{\left(\text{L}^\infty(\R)\right)^2}\right),
\end{equation}
where the constant $C$ has a dependence on $\delta$, $r$ and $R$ of the form $C=\cO(\delta^{-6}+R\delta^{-8})$
\end{rem}

\begin{coro}
Under the same assumptions as Theorem 1, the problem $(\widetilde{ \mathcal{H}})$ admits a unique solution $\wt u \in \text{H}^2_\loc(\wt \Omega)$, which can be approximated by $\widetilde{u}^{\text{app}}$ defined for almost every $(x,y)\in \wt \Omega$ by 
\begin{equation}
\widetilde{u}^{\text{app}}(x,y)=u^{\text{app}}\left(x,\frac{y}{h(x)}\right).
\end{equation}
Moreover, there exists a constant $\widetilde{C}>0$ depending only on $h_{\min}$, $h_{\max}$, $\delta$, $r$ and $R$ such that 
\begin{equation}
\Vert \wt u -\widetilde{u}^{\text{app}}\Vert_{\text{H}^1(\wt \Omega_r)}\leq \eta \wt C\left(\Vert \wt f\Vert_{\text{L}^2(\wt \Omega)}+\Vert \wt b \Vert_{(\text{H}^{1/2}(\partial\wt \Omega)}\right).
\end{equation}
\end{coro}

\begin{proof}
We use the equivalence between $(\widetilde{\mathcal{H}})$ and $(\mathcal{H})$, and we notice that 
\begin{equation}\nonumber
\Vert \wt u -\widetilde{u}^{\text{app}}\Vert_{\text{H}^1(\wt \Omega_r)}\leq h_{\max} \Vert u -u^{\text{app}}\Vert_{\text{H}^1(\Omega_r)}, \qquad \Vert \wt f\Vert_{\text{L}^2(\wt \Omega)}\geq h_{\min} \Vert f\Vert_{\text{L}^2(\Omega)}, 
\end{equation}
\begin{equation}\nonumber
\Vert \wt b \Vert_{(\text{H}^{1/2}(\R))^2}\geq h_{\min} \Vert b\Vert_{\left({\text{H}}^{1/2}(\partial\wt \Omega)\right)^2}.
\end{equation}
\end{proof}

\subsection{Modal Green functions and their approximations}

As mentioned in the previous section, we start by studying equations \eqref{multi1D} for every $n\in \N$. To this end, we denote by $G_n(x,s)$ the modal Green functions associated to $(\cH_n')$. It satisfies for every $s\in \R$ the partial differential equation
\begin{equation}
\left\{\begin{array}{cl}
\partial_{xx} G_n(x,s)+k_n(x)^2G_n(x,s)=-\delta_s & \text{ in } \R, \\ 
G_n(\cdot,s) \text{ is outgoing.} \end{array}\right. 
\end{equation}
We prove the following theorem which provides an approximation of $G_n$ for every $n\in N$ and the control of the approximation error in $\text{W}^{1,1}(\R)$. 

\begin{theorem}\label{th2}
For every $s\in \R$, the equation 
\begin{equation}\label{refgreen}
\left\{\begin{array}{lc}
\partial_{xx} G_n(x,s)+k_n(x)^2G_n(x,s)=-\delta_s & \text{ in } \R, \\ 
G_n(\cdot,s) \text{ is outgoing.} \end{array}\right. 
\end{equation}
has a unique solution $G_n(\cdot,s)\in \text{W}^{1,1}(\R)$. This solution can be decomposed as $G_n=G_n^{\text{app}}+\cO(\eta)$ where $G_n^{\text{app}}$ has the explicit form given in \eqref{greenfunction} and $\cO(\eta)$ is a term that tends to $0$ in $\text{W}^{1,1}(\R)$ uniformly in $s$ as $\eta$ tends to $0$. Moreover, let $r>0$. There exist $\eta_1>0$ depending on $R$, $r$, $h_{\min}$, $h_{\max}$ and $\delta$ such that if $\eta<\eta_1$, there exists $\alpha,\beta>0$ depending only on $h_{\min}$, $h_{\max}$, $r$, $\delta$ and $R$ such that for every $s\in \R$, 
\begin{equation}
\Vert G_n(\cdot,s) \Vert_{\text{L}^1(-r,r)}\leq \alpha_n^{(1)}:= \left\{\begin{array}{cl} \alpha & \text{ if } n\leq N, \\ \frac{\alpha}{\min(|k_n|)^2} & \text{ if } n> N, \end{array}\right. 
\end{equation} 
\begin{equation}
\Vert \partial_x G_n(\cdot,s) \Vert_{\text{L}^1(-r,r)}\leq \alpha_n^{(2)}:= \left\{\begin{array}{cl} \alpha & \text{ if } n\leq N, \\ \frac{\alpha}{\min(|k_n|)} & \text{ if } n> N, \end{array}\right. 
\end{equation} 
\begin{equation}
\Vert G_n(\cdot,s)-G_n^{\text{app}}(\cdot,s) \Vert_{\text{L}^1(-r,r)}\leq \beta_n^{(1)}:= \eta \left\{\begin{array}{cl} \beta & \text{ if } n\leq N,\\ \frac{\beta}{\min(|k_n|)^2} & \text{ if } n> N, \end{array}\right. 
\end{equation} 
\begin{equation}
\Vert \partial_x G_n(\cdot,s) -\partial_x G_n^{\text{app}}(\cdot,s)\Vert_{\text{L}^1(-r,r)}\leq \beta_n^{(2)} :=\eta \left\{\begin{array}{cl} \beta & \text{ if } n\leq N ,\\ \frac{\beta}{\min(|k_n|)} & \text{ if } n> N .\end{array}\right. 
\end{equation} 

\end{theorem}

\begin{rem}
The following proof shows that we can choose
\begin{equation}
\eta_1\leq \min\left(\frac{1}{9C_{N-1}},\frac{|k_{N+1}|}{9C_{N+1}},\frac{1}{4c_8\pi C_N}\right),
\end{equation}
where $C_{N-1}$ and $C_{N+1}$ are defined in \eqref{contrF}, $c_8$ is defined in Lemma \ref{lemconst} and $C_N$ comes from Theorem 2 in \cite{olver2}. It shows that $\eta$ has to be small compared to $r$, $R$, $h_{\min}$ and $\delta$ for this theorem to apply. 
\end{rem}

\begin{rem}
By looking at the proof, we can see that $\alpha$, $\beta$ and $C_n$ depend on $\delta$, $r$ and $R$ as $r\delta^{-1}$, $r R \delta^{-6}$, and $R\delta^{-5}$ respectively. 
\end{rem}

To prove this theorem, we first need a technical lemma to connect solutions of the partial differential equation \eqref{refgreen} defined for $x<s$ and $x>s$. 
\begin{lem}[Connection of the Green functions] \label{greenconnect}
Let $s\in \R$. Assume that $u$ is a solution to $u''+k_n^2 u=-\delta_s$ and that there exist $A,B\in \R$ and $w_1,w_2\in \mathcal{C}^2(\R)$ such that
\begin{equation}
u(x)=\left\{\begin{array}{cl} Aw_1(x) & \text{ if } x<s, \\ Bw_2(x) & \text{ if } x>s, \end{array}\right. 
\end{equation}
then 
\begin{equation}
A=\frac{w_2(s)}{w_1'(s)w_2(s)-w_2'(s)w_1(s)}, \qquad B=\frac{w_1(s)}{w_1'(s)w_2(s)-w_2'(s)w_1(s)}
\end{equation}
\end{lem}

\begin{proof}
Since $u$ is continuous in $s$, $Aw_1(s)=Bw_2(s)$. Then, using the jump formula for distributions, we find that 
\begin{equation}\nonumber
Bw_2'(s)-Aw_1'(s)=-1 \qquad \Rightarrow \qquad A(w_1(s)w_2'(s)-w_2(s)w_1'(s))=-w_2(s)
\end{equation}
\end{proof}

Next, we study the Green function for the three types of waves, depending on the value of $n$.

\subsection{Proof of Theorem 2}

\begin{proof}The propagative case ($n<N$)]

We denote $u_n=G_n(\cdot,s)$. Changing variable to $\sigma=\eta x$, we see that $w_n=u_n(\sigma/\eta)$ satisfies the equation $w_n''+k_n^2(\sigma/\eta)w_n^2/\eta^2=0$ for every $\sigma \neq \eta s$ where in the case at hand, $k_n^2(\sigma/\eta)>0$. Using Theorem 4 in \cite{olver1} on $w_n$, shows that there exist $A, B\in \C$ such that 
\begin{equation}\label{olv1}
u_n(x)=\left\{\begin{array}{cc}\frac{A}{\sqrt{k_n(x)}}\exp\left(- i\int_s^x k_n\right)(1+\overline{\eps(x)}) & \text{ if } x<s, \\ \frac{B}{\sqrt{k_n(x)}}\exp\left(i\int_s^x k_n\right)(1+\eps(x)) & \text{ if } x>s, \end{array}\right. 
\end{equation}
where $\eps\in \mathcal{C}^1(\R)$ is such that for all $x\in \R$,
\begin{equation}\label{eqeps}
|\eps(x)|\leq e^{F/2}-1, \quad |\eps'(x)|\leq 2k_n(x)(e^{F/2}-1),\end{equation}
and
\begin{equation} \nonumber |F|\leq \int_\R \eta \left| \frac{1}{\sqrt{k_n(\sigma/\eta)}}\partial_{\sigma\sigma} \left(\frac{1}{\sqrt{k_n(\sigma/\eta)}}\right)\right| \dd \sigma.
\end{equation}
Using the expression of $k_n$, we see that
\begin{equation}
\label{contrF}
|F| \leq \eta \int_{-R}^{R}\frac{1}{\eta^2} \left|\frac{n^2\pi^2(h''h-3(h')^2)}{2h^4k_n^3}+\frac{5(h')^2n^4\pi^4}{4h^6k_n^5}\right|(\sigma/\eta) \dd \sigma.
\end{equation}
We deduce that there exist a constant $\gamma_1>0$ depending on $h_{\min}$, $h_{\max}$ and $R$ such that
\begin{equation}\nonumber
|F|\leq \eta \gamma_1 \left(\frac{(N-1)^2\pi^2}{\delta^3}+\frac{(N-1)^4\pi^4}{\delta^5}\right):=\eta C_{N-1}.
\end{equation} 
Using Lemma \ref{greenconnect}, we find that 
\begin{equation}\nonumber
A=\frac{i(1+\eps(s))}{2\sqrt{k_n(s)}(1+R)}, \qquad B=\frac{i(1+\overline{\eps(s)})}{2\sqrt{k_n(s)}(1+R)},
\end{equation}
where 
\begin{equation}\nonumber R=2\text{Re}(\eps(s))+|\eps(s)|^2+\frac{\text{Im}(\eps'(s)(1+\overline{\eps(s)})}{k_n(s)}.
\end{equation}
It follows that 
\begin{equation}\nonumber
G_n(x,s)=\left\{\begin{array}{cc} G_n^{\text{app}}(x,s) \frac{(1+\overline{\eps(x)})(1+\eps(s))}{1+R} & \text{ if } x<s, \\
G_n^{\text{app}}(x,s) \frac{(1+\eps(x))(1+\overline{\eps(s)})}{1+R} & \text{ if } x>s, \end{array}\right.,\end{equation}
and 
\begin{multline}\nonumber
\partial_xG_n(x,s)= \partial_xG_n^{\text{app}}(x,s)\frac{1+\eps(s)}{1+R}\times \\
\left[\left(1+\overline{\eps(x)}+\frac{2\overline{\eps'(x)}k_n(x)^2h(x)^3}{2ik_n(x)^3h(x)^3-n^2\pi^2h'(x)}\right) \textbf{1}_{(-\infty,s)}(x)\right. \\ +\left. 
\left(1+\eps(x)+\frac{2\eps'(x)k_n(x)^2h(x)^3}{2ik_n(x)^3h(x)^3-n^2\pi^2h'(x)}\right)\textbf{1}_{(s,+\infty)}(x) \right].  
\end{multline}
Assuming that $\eta<\frac{1}{9C_{N-1}}$,then $|F|\leq 1/9$, 
\begin{equation}\nonumber
e^{F/2}-1\leq \frac{|F|/2}{1-|F|/4}\leq \min\left(\frac{3}{4}\eta C_n,\frac{1}{16}\right),
\end{equation}
and
\begin{equation}\nonumber
|R|\leq 2|\eps(x)|+|\eps|^2+\frac{|\eps'(x)|(1+|\eps(x)|)}{k_n(x)}\leq 4(e^{F/2}-1)e^{F/2}\leq \frac{1}{2}.
\end{equation}
It follows that 
\begin{equation}\nonumber
\left|1-\frac{(1+\overline{\eps(x)})(1+\eps(s))}{1+R}\right|\leq \frac{| \eps'(x)|(1+\Vert \eps\Vert_{L^\infty(\R)})/k_n(x)}{1-|R|}\leq 6\eta C_n,
\end{equation}
and that 
\begin{equation}\nonumber
\Vert G_n(\cdot,s)-G_n^{\text{app}}(\cdot,s)\Vert_{\text{L}^1(-r,r)}\leq \eta \frac{6C_n r }{\min(k_n)},\end{equation}
\begin{equation}\nonumber
\Vert G_n(\cdot,s) \Vert_{\text{L}^1(-r,r)}\leq \frac{2r}{\min(k_n)}.
\end{equation}
In the same way, 
\begin{equation}\nonumber
\left|\frac{1+\eps(s)}{1+R}\right|\left|\frac{2\eps'(x)k_n(x)^2h(x)^3}{2ik_n(x)^3h(x)^3-n^2\pi^2h'(x)}\right|\leq 4\frac{3}{2}\frac{C_n \eta k_n(x)^3 h(x)^3}{k_n(x)^3 h(x)^3}\leq 6C_n \eta
\end{equation}
and so 
\begin{multline} \nonumber\Vert \partial_xG_n(\cdot,s)-\partial_xG_n^{\text{app}}(\cdot,s)\Vert_{\text{L}^1(-r,r)} \leq \eta \frac{12C_n r }{\min(k_n)}\times \\ \left(\Vert k_n\Vert_{L^\infty(\R)}+\frac{n^2\pi^2}{18 C_n\min(k_n)^2h_{\min}^3}\right),
\end{multline}
and
\begin{equation}\nonumber
\Vert \partial_x G_n(\cdot,s) \Vert_{\text{L}^1(-r,r)}\leq \frac{3 r }{\min(k_n)}\left(\Vert k_n\Vert_{L^\infty(\R)}+\frac{n^2\pi^2}{18 C_n\min(k_n)^2h_{\min}^3}\right).
\end{equation}

\end{proof}

\begin{proof}[The evanescent case ($n>N$)]
We denote $u_n=G_n(\cdot,s)$. Changing variable to $\sigma=\eta x$, we see that $w_n=u_n(\sigma/\eta)$ satisfies the equation $w_n''+k_n^2(\sigma/\eta)w_n^2/\eta^2=0$ for every $\sigma \neq \eta s$. Theorem 3 in \cite{olver1} on $w_n$ yields the existence of $A, B\in \C$ such that 
\begin{equation}\label{olv2}
u_n(x)=\left\{\begin{array}{cc} \frac{A}{\sqrt{|k_n|(x)}}\exp\left(\int_s^x |k_n|\right)(1+\eps_2(x)) & \text{ if } x<s, \\ \frac{B}{\sqrt{|k_n|(x)}}\exp\left(-\int_s^x |k_n|\right)(1+\eps_1(x)) & \text{ if } x>s,\end{array}\right. 
\end{equation}
where $\eps_1,\eps_2\in \mathcal{C}^1(\R)$ are such that for $i=1,2$ and all $x\in \R$,
\begin{equation}\label{eqeps2}
|\eps_i(x)|\leq e^{F/2}-1, \quad |\eps'_i(x)|\leq 2|k_n|(x)(e^{F/2}-1),
\end{equation}
where $F$ satisfies \eqref{contrF}. It follows that there exist a constant $\gamma_2>0$ depending on $h_{\min}$, $h_{\max}$ and $R$ such that 
\begin{equation}\nonumber
|F|\leq \frac{\eta \gamma_2 (N+1)^2\pi^2}{\min(|k_n|)\delta^2}\left(1+\frac{(N+1)^2\pi^2}{\delta^2}\right):=\eta \frac{C_{N+1}}{\min(|k_n|)}.
\end{equation}
Using Lemma \ref{greenconnect}, we find that 
\begin{equation}\nonumber
A=\frac{1+\eps_1(s)}{2\sqrt{k_n(s)}(1+R)}, \qquad B=\frac{1+\eps_2(s)}{2\sqrt{k_n(s)}(1+R)},
\end{equation}
where 
\begin{equation}\nonumber R=\left(\eps_2+\eps_1+\eps_1\eps_2+\frac{\eps_2'(1+\eps_1)}{2|k_n|}-\frac{\eps_1'(1+\eps_2)}{2|k_n|}\right)(s).
\end{equation}
It follows that 
\begin{equation}\nonumber
G_n(x,s)=\left\{\begin{array}{cc} G_n^{\text{app}}(x,s) \frac{(1+\eps_2(x))(1+\eps_1(s))}{1+R} & \text{ if } x<s, \\
G_n^{\text{app}}(x,s) \frac{(1+\eps_1(x))(1+\eps_2(s))}{1+R} & \text{ if } x>s, \end{array}\right.\end{equation}
and 
\begin{multline}\nonumber
\partial_x G_n(x,s)= \partial_x G_n^{\text{app}}(x,s)\times \\
\left[\frac{1+\eps_1(s)}{1+R}\left(1+\eps_2(x)-\frac{2\eps'_2(x)|k_n(x)|^2h(x)^3}{h'(x)n^2\pi^2+2|k_n(x)|^3h(x)^3}\right) \textbf{1}_{(-\infty,s)}(x) \right. \\
+ \left. \frac{1+\eps_2(s)}{1+R}\left(1+\eps_1(x)-\frac{2\eps'_1(x)|k_n(x)|^2h(x)^3}{h'(x)n^3\pi^2+2|k_n(x)|^3h(x)^3}\right) \textbf{1}_{(s,+\infty)}(x)\right] . 
\end{multline}
We also notice that 
\begin{equation}\nonumber
\left\Vert e^{- \int_s^x |k_n|}\right\Vert_{\text{L}^1(-r,r)}\leq \Vert e^{- \min(|k_n|)|x-s| }\Vert_{\text{L}^1(-r,r)}\leq \frac{2}{\min(|k_n|)}.\end{equation}
Assuming that $\eta<\frac{|k_{N+1}|}{9C_{N+1}}$ then $|F|\leq 1/9$, 
\begin{equation}\nonumber
e^{F/2}-1\leq \frac{|F|/2}{1-|F|/4}\leq \min\left(\frac{3}{4}\frac{\eta C_{N+1}}{\min(|k_{n}|)},\frac{1}{16}\right),
\end{equation}
and
\begin{equation}\nonumber
|R|\leq 2(e^{F/2}-1)+(e^{F/2}-1)^2+2e^{F/2}(e^{F/2}-1)\leq 4(e^{F/2}-1)e^{F/2}\leq \frac{1}{2}.
\end{equation}
It follows that 
\begin{multline}\nonumber
\left|1-\frac{(1+\eps_1(x))(1+\eps_2(s))}{1+R}\right|\leq \\ \frac{| \eps_1(x)|(1+\Vert \eps_2\Vert_{L^\infty(\R)})+| \eps_2(x)|(1+\Vert \eps_1\Vert_{L^\infty(\R)})}{2|k_n(x)|(1-|R|)}\leq \frac{6\eta C_{N+1}}{\min(|k_n|)}.
\end{multline}
and that \begin{equation}\label{c1}
\Vert G_n(\cdot,s)-G_n^{\text{app}}(\cdot,s)\Vert_{\text{L}^1(-r,r)}\leq \eta\frac{12C_{N+1}r}{\min(|k_n|)^2},\end{equation}
\begin{equation}\label{c2}
\Vert G_n(\cdot,s) \Vert_{\text{L}^1(-r,r)}\leq  \frac{4r}{\min(|k_n|)^2},
\end{equation}
In the same way, 
\begin{multline}\nonumber
\left|\frac{1+\eps_1(s)}{1+R}\right|\left|\frac{2\eps_2'(x)|k_n(x)|^2h(x)^3}{h'(x)n^2\pi^2+2|k_n(x)|^3h(x)^3}\right| \\ \leq 4\frac{3}{2}\frac{2C_{N+1} \eta |k_n(x)|^3h(x)^3/\min(|k_n|)}{2|k_n(x)|^3h(x)^3}\leq \frac{6C_{N+1} \eta }{\min(|k_n|)},
\end{multline}
and so 
\begin{equation}\label{c3}\Vert \partial_x G_n(\cdot,s)-\partial_xG_n^{\text{app}}(\cdot,s)\Vert_{\text{L}^1(-r,r)} \leq \eta \frac{24r}{\min(|k_n|)^2}\left(\Vert k_n\Vert_{L^\infty(\R)} +\frac{(N+1)^2\pi^2}{2\delta^2h_{\min}^3}\right),
\end{equation}
and 
\begin{equation}\label{c4}
\Vert \partial_x G_n(\cdot,s) \Vert_{\text{L}^1(-r,r)}\leq \frac{6r}{\min(|k_n|)^2}\left(\Vert k_n\Vert_{L^\infty(\R)} +\frac{(N+1)^2\pi^2}{2\delta^2h_{\min}^3}\right).
\end{equation}
\end{proof}

\begin{rem}
We notice that the control of $G_n$, $\partial_x G_n$ and $\partial_x G_n -\partial_xG_n^{\text{app}}$ is uniform in $n>N$, with
\begin{equation}
\eta<\frac{|k_{N+1}|}{9C_{N+1}}.
\end{equation}
This uniform control is essential to obtain the global control $\eta<\eta_1$ in Theorem \ref{th2}. We could also provide a uniform control in inequalities \eqref{c1}, \eqref{c2}, \eqref{c3} and \eqref{c4}. However, in the following, we add these inequalities and thus we keep track of the factors $1/\min(|k_n|)$ to ensure fast decrease when $n$ goes to infinity. 
\end{rem}

It remains to deal with the case $n=N$. This case is more complicated, since $x\mapsto k_n(x)^2$ is not of constant sign. We first introduce two technical Lemmas used to give an approximation of the Green function for $n=N$.

\begin{lem}\label{lemxi}
Let us define
\begin{equation}
\xi(x)=\left\{\begin{array}{cc} \left(-\frac{3}{2}i\int_x^{x^*}k_N\right)^{2/3} & \text{ if } x<x^*, \\ -\left(\frac{3}{2}\int_{x^*}^x k_N\right)^{2/3} & \text{ if } x>x^*.\end{array}\right. 
\end{equation}
This function is a decreasing bijection from $\R$ to $\R$. Moreover, the function 
\begin{equation} \phi : x\mapsto (-\xi(x))^{1/4}/\sqrt{k_n(x)},\end{equation} is in $\mathcal{C}^2(\R)$ and for all $x\in \R$, there exists a constant $c_\phi$ such that $|\phi'(x)|\leq c_\phi|\phi(x)|$.
\end{lem}

\begin{proof}
This lemma can be proved by adapting section 4 of \cite{olver2}. 
\end{proof}

\begin{lem}\label{lemconst}
Let us define the following functions:
\begin{equation}
E(x)=\exp\left(\frac{2}{3}x^{3/2}\right)1_{x>0} +1_{x\leq 0},
\end{equation}
\begin{equation}
M=\sqrt{E^2\Ai^2+E^{-2}\Bi^2}, \qquad N=\sqrt{E^2(\Ai')^2+E^{-2}(\Bi')^2}
\end{equation}
There exist a constant $c_1\approx 0.4$ such that for all $x\in \R$,
\begin{equation}
|M(x)E(x)\Ai'(x)|\leq c_1, \qquad \left|\frac{\Bi(x)N(x)}{E(x)}\right|\leq c_1, \qquad |M(x)N(x)|\leq c_1,
\end{equation}
\begin{equation}
\left|\frac{\Bi'(x)M(x)}{E(x)}\right|\leq c_1, \qquad |\Ai(x)E(x)N(x)|\leq c_1.
\end{equation}
There also exist constants $c_2, c_3$ depending on $h_{\min}$, $h_{\max}$, $k$ and $R$ such that for every $x\in \R$,
\begin{equation}
\left|\frac{(-\xi(x))^{1/4}}{\sqrt{k_n(x)}}M(\xi(x))\right|\leq c_6, \qquad \left|\frac{\sqrt{k_n(x)}}{(-\xi(x))^{1/4}}N(\xi(x))\right|\leq c_7.
\end{equation}

\end{lem}
\begin{proof}
Using Airy's asymptotic expansions presented in section 10.4 of \cite{abramowitz1}, we obtain the first constant and the control
\begin{equation}\nonumber
M(x)=\cO_{|x|\rightarrow \infty}\left(\frac{1}{|x|^{1/4}}\right), \qquad N(x)=\cO_{|x|\rightarrow \infty}\left(|x|^{1/4}\right). 
\end{equation}
Since $h'$ is compactly supported, we conclude the proof by noticing that 
\begin{equation}\nonumber
k_n(x)=\cO_{|x|\rightarrow \infty}(1), \qquad \xi(x)=\cO_{|x|\rightarrow \infty}(|x|).
\end{equation}

\end{proof}

\begin{proof}[The locally resonant case ($n=N$)]
We set $u_n=G_n(\cdot,s)$. Changing variable to $\sigma=\eta x$, we see that $w_n=u_n(\sigma/\eta)$ satisfies the equation $w_n''+k_n^2(\sigma/\eta)w_n^2/\eta^2=0$ for every $\sigma \neq \eta s$. This equation is very similar to the Airy equation and from Theorem 2 in \cite{olver2}, we know that there exist $A, B\in \C$ such that 
\begin{equation}\label{olv3}
u_n(x)=\left\{\begin{array}{cc} \frac{A(-\xi(x))^{1/4}}{\sqrt{k_n(x)}}(\Ai+\eps_1)(\xi(x)) & \text{ if } x<s, \\ \frac{B(-\xi(x))^{1/4}}{\sqrt{k_n(x)}}(i\Ai+\Bi+i\eps_1+\eps_2)(\xi(x)) & \text{ if } x>s,\end{array}\right. 
\end{equation}
where $\eps_1,\eps_2\in \mathcal{C}^1(\R)$ are such that there exist $C_N$ depending on $h_{\min}$, $h_{\max}$, $\delta$ and $R$ such that for all $x\in \R$,
\begin{equation}\label{eqeps3}
\left|\frac{E(\xi)}{M(\xi)}\right||\eps_1(\xi)|\leq \frac{1}{\lambda_1}(e^{\eta C_N\lambda_1}-1), \quad \left|\frac{E(\xi)}{N(\xi)}\right||\eps_1'(\xi)|\leq \frac{1}{\lambda_1} (e^{\eta C_N\lambda_1}-1),
\end{equation}
\begin{equation}\label{eqeps4}
\left|\frac{1}{E(\xi)M(\xi)}\right||\eps_2(\xi)|\leq \frac{\lambda_2}{\lambda_1}(e^{\eta C_N\lambda_1}-1), \quad \left|\frac{1}{E(\xi)N(\xi)}\right||\eps_2'(\xi)|\leq \frac{\lambda_2}{\lambda_1}(e^{\eta C_N\lambda_1}-1),
\end{equation}
where $\lambda_1>\lambda_2$ are known constants. Using Lemma \ref{greenconnect} and the fact that $\Bi\Ai'-\Bi'\Ai=-1/\pi$, we find that 
\begin{equation}\nonumber
A=\frac{-\pi(-\xi(s))^{1/4}(i\Ai(\xi(s))+\Bi(\xi(s))+i\eps_1(\xi(s))+\eps_2(\xi(s)))}{(1-R\pi)\sqrt{k_n(s)}},\end{equation}
\begin{equation}\nonumber
B=\frac{-\pi(-\xi(s))^{1/4}(\Ai(\xi(s))+\eps_1(\xi(s)))}{(1-R\pi)\sqrt{k_n(s)}},
\end{equation}
where 
\begin{equation}\nonumber
R=(\Ai'\eps_2+\eps_1'\Bi+\eps_1'\eps_2-\Bi'\eps_1-\Ai\eps_2'-\eps_2'\eps_1)(\xi(s)).
\end{equation}
It follows that 
\begin{multline}\nonumber
G_n(x,s)= G_n^{\text{app}}(x,s)\times \\
\left[\left(1+\frac{\eps_1(\xi(x))}{\Ai(\xi(x))}\right)\frac{1}{1-R\pi}\left(1+\frac{i\eps_1(\xi(s))+\eps_2(\xi(s))}{i\Ai(\xi(s))+\Bi(\xi(s))}\right)\textbf{1}_{(-\infty,s)}(x)\right.  \\
 +\left. \left(1+\frac{i\eps_1(\xi(x))+\eps_2(\xi(x))}{i\Ai(\xi(x))+\Bi(\xi(x))}\right)\frac{1}{1-R\pi}\left(1+\frac{\eps_1(\xi(s))}{\Ai(\xi(s))}\right) \textbf{1}_{(s,+\infty)}(x)\right],
\end{multline}
and 
\begin{multline}\nonumber
\partial_x G_n^{\text{app}}(x,s)=\frac{-\pi(-\xi(s))^{1/4}}{(1-R\pi)\sqrt{k_n(s)}}\times \\
 \left[\left((\Ai'+\eps_1')(\xi(x))\frac{\sqrt{k_n(x)}}{(-\xi(x))^{1/4}}+(\Ai+\eps_1)(\xi(x))\phi'(x)\right)  \right. \\
 \times (i\Ai+\Bi+i\eps_1+\eps_2)(\xi(s)) \textbf{1}_{(-\infty,s)}(x) \\
 + \left((i\Ai'+\Bi+i\eps_1'+\eps_2')(\xi(x))\frac{\sqrt{k_n(x)}}{(-\xi(x))^{1/4}}+(i\Ai
+\Bi+i\eps_1+\eps_2)(\xi(x))\phi'(x)\right) \\
 \times (\Ai+\eps_1)(\xi(s)) \textbf{1}_{(x,+\infty)}(x)\Bigg].
 \end{multline}
We define $c_4=(4\lambda_2+2)c_1$. Assuming that $\eta\leq 1/(4C_Nc_4\pi)$, we know that $\eta\leq 1/(C_N\lambda_1)$ and $(e^{\eta C_N \lambda_1}-1)/\lambda_1\leq 2\eta C_N$, so 
\begin{equation}\nonumber
|R|\leq 2\eta C_N c_4, \qquad \left|\frac{R\pi}{1-R\pi}\right|\leq \frac{2\eta C_N c_4\pi }{1-2\eta \pi C_N c_4} \leq 4\eta \pi  C_N c_4\leq 1,
\end{equation}
\begin{multline}\nonumber
\left|\frac{(-\xi(s))^{1/4}(-\xi(x))^{1/4}}{\sqrt{k_n(s)}\sqrt{k_n(x)}}\Ai(\xi(x))(i\eps_1(\xi(s))+\eps_2(\xi(s)))\right|\leq\qquad \qquad \qquad  \\ 
\left|\frac{(-\xi(s))^{1/4}(-\xi(x))^{1/4}}{\sqrt{k_n(s)}\sqrt{k_n(x)}}\right|4\lambda_2\eta C_N M(\xi(x))M(\xi(s))\leq 4\lambda_2\eta C_N c_2^2,
\end{multline}
\begin{multline}\nonumber
\left|\frac{(-\xi(s))^{1/4}(-\xi(x))^{1/4}}{\sqrt{k_n(s)}\sqrt{k_n(x)}}\eps_1(\xi(x))\left(i\Ai(\xi(s))+\Bi(\xi(s))+i\eps_1(\xi(s))+\eps_2(\xi(s))\right)\right|\\ \leq (1+2\lambda_2\eta C_N)4\eta C_N c_2^2.
\end{multline}
It follows that
\begin{equation}\nonumber
|G_n(x,s)-G_n^{\text{app}}(x,s)|\leq 4\eta \pi C_N c_4 |G_n^{\text{app}}(x)|+ \left(8\pi c_2^2\lambda_2\eta C_N + 24\pi \eta C_N c_2^2\right),
\end{equation}
which leads to 
\begin{equation}\nonumber
|G_n(x)-G_n^{\text{app}}(x)|\leq 4 \eta C_N \pi \left( c_4 |G_n^{\text{app}}(x)|+c_2^2\left(2\lambda_2 + 6\right)\right).
\end{equation}
Using the same idea, we prove that 
\begin{equation}\nonumber
|G_n^{\text{app}}(x,s)|\leq 2c_2^2\pi,
\end{equation}
and it follows that
\begin{equation}\nonumber
\Vert G_n(\cdot,s)-G_n^{\text{app}}(\cdot,s) \Vert_{\text{L}^1(-r,r)}\leq \eta 16r C_N\pi c_2^2 (c_4\pi+\lambda_2+3),\end{equation}
and
\begin{equation}\nonumber \Vert G_n(\cdot,s)\Vert_{\text{L}^1(-r,r)}\leq 4rc_2^2\pi\left(2+\frac{\lambda_2+3}{c_4\pi}\right).
\end{equation}
Using the same technique also prove that 
\begin{equation}\nonumber
\Vert \partial_xG_n(\cdot,s)-\partial_xG_n^{\text{app}}(\cdot,s)\Vert_{\text{L}^1(-r,r)}\leq \eta 8r C_N\pi c_2(c_3+c_2c_\phi)(2c_4+3\lambda_2),
\end{equation}
\begin{equation}\nonumber
\Vert \partial_x G_n(\cdot,s)\Vert_{\text{L}^1(-r,r)}\leq 16r\pi c_2(c_3+c_2c_\phi).
\end{equation}
\end{proof}

\subsection{Proof of Theorem \ref{th1}}

As mentioned at the beginning of section 3.1, we map the deformed waveguide $\widetilde{\Omega}$ to the regular waveguide $\Omega$ by a change of variables $\psi: (x,y)\mapsto (x,h(x) y)$. The problem \eqref{MAFT} is equivalent in $\Omega$ to the problem 
\begin{equation} \label{eqomega} \left\{\begin{array}{cl}
\partial_{xx} u +k^2 u +\frac{1}{h^2}\partial_{yy} u -\frac{h''h-2(h')^2}{h^3}y\partial_y u \hspace{1.5cm}& \\
\hspace{1.5cm}+\frac{(h')^2}{h^4}y^2\partial_{yy}u -\frac{2h'}{h^2}y \partial_{yx} u =-f &\text{ in } \Omega, \\ \partial_\nu u =b_\top+\frac{h'}{h}\partial_x u & \text{ on }\partial\Omega_{\text{top}}, \\ \partial_\nu u =b_\bot & \text{ on } \partial\Omega_{\text{bot}}, \\ u \text{ is outgoing.} \end{array}\right. \tag{$\mathcal{H}$}
\end{equation}
If we try and use the modal decomposition on this equation, mode coupling appears. We can however try to approach the solutions of $(\cH)$ by the solutions of the following system 
\begin{equation} 
\left\{\begin{array}{cl}\displaystyle \partial_{xx} v +\frac{1}{h(x)^2}\partial_{yy} v+k^2 v = -f & \text{ in } \Omega, \\ \partial_\nu v = b_\top & \text{ on } \partial\Omega_\text{top},\\ \partial_\nu v = b_\bot & \text{ on } \partial\Omega_\text{bot}, \\ v \text{ is outgoing}, \end{array}\right. \tag{$\mathcal{H}'$}
\end{equation}
which is amenable to modal decomposition. To estimate the error of such an approximation, we need to control the dependence between the source and the solution of $(\cH')$. In Proposition \ref{contrs}, we provide a control of the wave field generated by a source term in the waveguide, and in Proposition \ref{contrb}, we do the same thing for a source term on its boundary. The proofs of both propositions are given in the appendix.
\begin{rem} $(\cH')$ was obtained from $(\cH)$ by formally eliminating the terms likely to cause mode coupling. We cannot however neglect the term $\frac{1}{h(x)^2}\partial_{yy} u$ and approximate 
\begin{equation}\nonumber
\frac{1}{h(x)^2}\approx \frac{1}{h(x_0)^2}+\eta \,\cO_{x\to x_0}(x-x_0),
\end{equation}
for a constant $x_0\in \R$. Indeed, if $x$ is large enough, $x-x_0\gg \eta$ and 
\begin{equation}\nonumber
\frac{1}{h(x_0)^2}-\frac{1}{h_{\min}^2}=\cO(1), \qquad \frac{1}{h(x_0)^2}-\frac{1}{h_{\max}^2}=\cO(1).
\end{equation}
\end{rem}

\begin{prop}\label{contrs}
Let $r>0$ and $f\in \text{L}^2(\Omega_r)$. The equation
\begin{equation}\label{hdroit}
\left\{\begin{array}{cc} \displaystyle \partial_{xx} u+\frac{1}{h^2}\partial_{yy} u+k^2u=-f & \text{ in } \Omega, \\ \partial_\nu u =0 & \text{ on } \partial\Omega, \\ u \text{ is outgoing,} \end{array}\right.
\end{equation}
has a unique solution $u\in \text{H}^2_\loc(\Omega)$. Using notations of Theorem \ref{th2}, if $\eta<\eta_1$ then the operator 
\begin{equation}
\Gamma: \begin{array}{rcl} \text{L}^2(\Omega_r) & \rightarrow & \text{H}^2(\Omega_r) \\ f & \mapsto & u_{|\Omega_r} \end{array}, \quad \text{ where } u \text{ is the  solution to } \eqref{hdroit},
\end{equation}
is well defined, continuous and there exists a constant $D_1$ depending on $\delta$, $h_{\min}$, $h_{\max}$, $R$ and $r$ such that
\begin{equation}
\Vert u \Vert_{\text{H}^2(\Omega_r)}\leq D_1\Vert f \Vert_{\text{L}^2(\Omega_r)}.
\end{equation}
\end{prop}

\begin{prop}\label{contrb}
Let $r>0$ and $b=(b_\top,b_\bot)\in (\widetilde{\text{H}}^{1/2}(-r,r))^2$. The equation
\begin{equation}\label{hdroitb}
\left\{\begin{array}{cc} \displaystyle \partial_{xx} u+\frac{1}{h^2}\partial_{yy} u+k^2u=0 & \text{ in } \Omega, \\ \partial_\nu u =b_\top & \text{ on } \partial\Omega_\text{top}, \\ \partial_\nu u= b_\bot & \text{ on } \partial\Omega_\text{bot}, \\ u \text{ is outgoing.} \end{array}\right.
\end{equation}
has a unique solution $u\in \text{H}^2_\loc(\Omega)$. Using notations of Theorem \ref{th2}, if $\eta<\eta_1$ then 
the operator 
\begin{equation}
\Pi: \begin{array}{rcl} \left(\widetilde{\text{H}}^{1/2}(-r,r)\right)^2 & \rightarrow & \text{H}^2(\Omega_r) \\ b=(b_\top,b_\bot) & \mapsto & u_{|\Omega_r} \end{array}, \quad \text{ where } u \text{ is the solution to } \eqref{hdroitb},
\end{equation}
is well defined, continuous and there exists a constant $D_2$ depending on $\delta$, $h_{\min}$, $h_{max}$, $R$ and $r$ such that 
\begin{equation}
\Vert u \Vert_{\text{H}^2(\Omega_r)}\leq D_2\Vert b \Vert_{\widetilde{(\text{H}}^{1/2}(-r,r))^2}.
\end{equation}
\end{prop}

Using these two propositions, we are now able to justify the approximation of $(\cH)$ by $(\cH')$ which, as in \cite{colton1} and \cite{bonnetier1}, is a Born approximation. However, here we show that this approximation remains valid near resonance frequencies. 

\begin{prop}[Born approximation] \label{propborn}
Let $\mathcal{S}: \text{H}^2(\Omega_r)\rightarrow \text{L}^2(\Omega_r)$ and $\mathcal{T}: \text{H}^2(\Omega_r)\rightarrow (\widetilde{\text{H}}^{1/2}(-r,r))^2$ and $D_1$, $D_2$ the constants defined in Propositions \ref{contrs} and \ref{contrb}. Let $f\in \text{L}^2(\Omega_r)$ and $b\in (\widetilde{\text{H}}^{1/2}(-r,r))^2$. If
\begin{equation}\label{bornhyp}
\mu:= D_1\Vert \mathcal{S}\Vert_{\text{H}^2(\Omega_r),\text{L}^2(\Omega_r)}+D_2\Vert \mathcal{T}\Vert_{\text{H}^2(\Omega_r), (\widetilde{\text{H}}^{1/2}(-r,r))^2}<1,
\end{equation}
then the equation 
\begin{equation}\label{eqborn}
u=\Gamma(f)+\Pi(b)+\Gamma(\mathcal{S}(u))+\Pi(\mathcal{T}(u)),\end{equation}
has a unique solution $u\in \text{H}^2(\Omega_r)$. Moreover, if we define $v=\Gamma(f)+\Pi(b)$ then 
\begin{equation}
\Vert u-v\Vert_{\text{H}^2(\Omega_r)}\leq  \left(D_1\Vert f\Vert_{\text{L}^2(\Omega_r)}+D_2\Vert b\Vert_{(\widetilde{\text{H}}^{1/2}(-r,r))^2}\right)\frac{\mu}{1-\mu}.
\end{equation}
\end{prop}
\begin{proof}
If \eqref{bornhyp} is satisfied then $\Gamma\circ \mathcal{S}+\Pi\circ \mathcal{T}$ is a contraction and $u$ can be expressed into a Born series (see \cite{colton1}). We conclude using the results on geometrical series. 
\end{proof}

Coming back to equation \eqref{eqomega}, we define the operators 
\begin{equation}\nonumber
\mathcal{S} : \begin{array}{rcl} \text{H}^2(\Omega_r) & \rightarrow & \text{L}^2(\Omega_r) \\
u & \mapsto & \frac{h''h-2(h')^2}{h^3}y\partial_y u -\frac{(h')^2}{h^4}y^2\partial_{yy}u +\frac{2h'}{h^2}y \partial_{yx} u \end{array} ,
\end{equation}
and 
\begin{equation}\nonumber
\mathcal{T} : \begin{array}{rcl} \text{H}^2(\Omega_r) & \rightarrow & \widetilde{\text{H}}^{1/2}(-r,r) \\
u & \mapsto & \frac{h'}{h}\partial_xu|_{y=1} \end{array} .
\end{equation}
With these definitions, \eqref{eqomega} can be rewritten as \eqref{eqborn}. We define  
\begin{equation}\nonumber
\eta_0= \min\left(1,\eta_1,\frac{1}{2 \frac{D_1}{{h_{\min}}^2} \left(3+\frac{1}{{h_{\min}}^2}+\frac{2}{h_{\min}}\right)+2D_2\frac{1}{h_{\min}}}\right),
\end{equation}
and we notice that if $\eta \leq \eta_0$ then 
\begin{equation}\nonumber
\Vert \mathcal{S}\Vert_{\text{H}^2(\Omega_r), \text{L}^2(\Omega_r)}\leq \frac{\eta}{{h_{\min}}^2} \left(3+\frac{1}{{h_{\min}}^2}+\frac{2}{h_{\min}}\right),
\end{equation}
\begin{equation}\nonumber
\Vert \mathcal{T}\Vert_{\text{H}^2(\Omega_r), (\widetilde{\text{H}}^{1/2}(-r,r))^2}\leq \eta \frac{1}{h_{\min}}. 
\end{equation}
We choose $\mu$ as in \eqref{bornhyp} and we define $v=\Gamma(f)+\Pi(b_\top,b_\bot)$. Using Proposition \ref{propborn}, the problem \eqref{eqomega} has a unique solution $u\in \text{H}^2(\Omega_r)$ and 
\begin{multline}\nonumber
\Vert u-v\Vert_{\text{H}^2(-r,r)}\leq 2\eta \left(\frac{D_1}{{h_{\min}}^2} \left(3+\frac{1}{{h_{\min}}^2}+\frac{2}{h_{\min}}\right)+D_2\frac{1}{h_{\min}}\right) \\ \times\left(D_1\Vert f\Vert_{\text{L}^2(\Omega_r)}+D_2\Vert b\Vert_{(\widetilde{\text{H}}^{1/2}(-r,r))^2}\right).
\end{multline}
Using the modal decomposition, we know that 
\begin{equation}\nonumber
v(x,y)=\sum_{n\in \N} \left(\int_\R G_n(x,s) (f_n(s)+b_\top(s)\ph_n(1)+b_\bot(s)\ph_n(0))\dd s\right) \ph_n(y).
\end{equation}

We now estimate the error between $v$ and $u^{\text{app}}$ given by \eqref{greentot}, following the same idea as in the proof of Proposition \ref{contrs}. We denote $g_n=f_n+b_\top\ph_n(1)+b_\bot\ph_n(0)$ and $g=\tsum_{n\in \N} g_n\ph_n$. Using the Young inequality for integral operators and the results and notations of Theorem \ref{th2},
\begin{equation}\nonumber
\Vert v_n-u^{\text{app}}_n\Vert_{\text{L}^2(-r,r)}\leq \beta_n^{(1)}\Vert g_n\Vert_{\text{L}^2(-r,r)}, \,\,\, \Vert v_n'-(u^{\text{app}}_n)'\Vert_{\text{L}^2(-r,r)}\leq \beta_n^{(2)}\Vert g_n\Vert_{\text{L}^2(-r,r)}.
\end{equation}
It follows that 
\begin{multline}\nonumber
\Vert v-u^{\text{app}}\Vert_{\text{H}^1(\Omega_r)}^2\leq \sum_{n=0}^{N} (2+n^2\pi^2)\beta^2\eta^2\Vert g_n\Vert^2_{\text{L}^2(-r,r)}\\+\sum_{n>N} \left(\frac{(1+n^2\pi^2)\beta^2\eta ^2}{\min(|k_{n}|)^4}+\frac{\eta^2\beta^2}{\min(|k_{n}|)^2}\right)\Vert g_n\Vert_{\text{L}^2(-r,r)}^2,
\end{multline}
and so
\begin{equation}\nonumber
\Vert v-u^{\text{app}}\Vert_{\text{H}^1(\Omega_r)}^2\leq \beta^2 \eta^2 \max\left(2+N^2\pi^2,\frac{1+(N+1)^2\pi^2}{\delta^4}+\frac{1}{\delta^2}\right) \Vert g\Vert_{\text{L}^2(\Omega_r)}^2.
\end{equation}
We conclude by noticing that 
\begin{equation}\nonumber
\Vert u-u^{\text{app}}\Vert_{\text{H}^1(\Omega_r)}\leq \Vert u-v\Vert_{\text{H}^2(\Omega_r)}+\Vert v-u^{\text{app}}\Vert_{\text{H}^1(\Omega_r)}.
\end{equation}

\section{Extension to general slowly varying waveguides \label{section4}}

\subsection{The cut and match strategy}
In the previous section, we constructed an approximation for the solution to the Helmholtz equation in a slowly increasing waveguide. In this section, we generalize our result by considering a 2D infinite waveguide $\widetilde{\Omega}=\{(x,y)\in \R^2 \,|\, 0< y < h(x)\}$ where $h\in \mathcal{C}^2(\R)$ is such that $h'$ is compactly supported and there exists a parameter $\eta$, assumed to be small compared to $R$ and $h_{\min}$, such that $\Vert h'\Vert_{L^\infty(\R)} \leq \eta$ and $\Vert h'' \Vert_{L^\infty(\R)} \leq \eta^2$. 

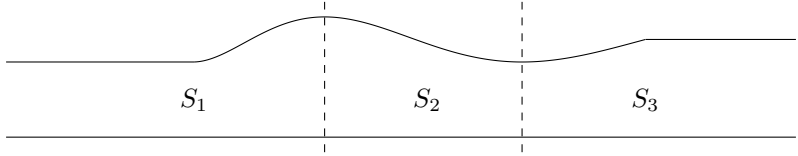
\begin{figure}[h]
\begin{center}
\begin{tikzpicture} 
\draw (-1.5,0) -- (9,0);
\draw (-1.5,1) -- (1,1); 
\draw (7,1.2998) -- (9,1.2998);  
\draw [domain=1:7, samples=200] plot (\x,{0.3*sin(4*pi*sqrt(\x)/sqrt(7) r)+1-0.3*sin(4*pi/sqrt(7) r)});  
\draw [dashed] (2.7344, -0.2)--(2.7344,1.8); 
\draw [dashed] (5.3594,-0.2)--(5.3594,1.8); 
\draw (1,0.5) node{$S_1$}; 
\draw (4.1,0.5) node{$S_2$}; 
\draw (7,0.5) node{$S_3$}; 
\end{tikzpicture}
\end{center}
\caption{\label{fig2} Representation of the sections in a slowly varying waveguide.}
\end{figure}

When we look at the proof of Theorem \ref{th1}, we notice that the condition that $h$ is increasing is only required to properly define the change of variable $x\mapsto  \xi(x)$ in Theorem \ref{th2} when $n=N$. In order to generalize Theorem \ref{th1}, we only need to generalize Theorem \ref{th2} for the case $n=N$. To this end, we follow a strategy developed in the context of the Schrödinger equation, see for instance \cite{roy1}. We partition $\widetilde{\Omega}$ in $J$ regions $(S_j)_{j=1,..,J}$,  on which $h'$ has a constant sign, as shown in Figure \ref{fig2}.

If there exists $x\in \R$ such that $k_n(x)=0$ in $S_j$, we denote this coordinate by $x^*_j$. Otherwise, as in Section 5 of \cite{olver2}, $x^*_j$ is chosen to be greater than $\max(S_j)$ if $k_n^2$ is positive on $S_j$ (resp. smaller than $\min(S_j)$ if $k_n^2$ is negative on $S_j$). Then, if $h$ is increasing in $S_j$, we define
\begin{equation}
\xi_j(x)=\left\{\begin{array}{cc} \left(-\frac{3}{2}i\int_x^{x^*_j}k_n(t)\dd t\right)^{2/3} & \text{ if } x<x^*_j \text{ and } x\in S_j, \\  -\left(\frac{3}{2}\int_{x^*_j}^x k_n(t) \dd t \right)^{2/3} & \text{ if } x>x^*_j \text{ and } x\in S_j.\end{array}\right.
\end{equation}
Otherwise, if $h$ is decreasing in $S_j$, we define
\begin{equation}
\xi_{j}(x)=\left\{\begin{array}{cc} -\left(\frac{3}{2}\int_x^{x^*_j}k_n(t)\dd t\right)^{2/3} & \text{ if } x<x^*_{j} \text{ and } x\in S_j, \\  \left(-\frac{3}{2}i\int_{x^*_{j}}^x k_n(t) \dd t \right)^{2/3} & \text{ if } x>x^*_{j} \text{ and } x\in S_j.\end{array}\right.
\end{equation}
In both cases, we denote
\begin{equation}\label{phii}
\phi_j(x)=\frac{(-\xi_j(x))^{1/4}}{\sqrt{k_n(x)}}.
\end{equation}

Given $s\in \R$, we study the problem 
\begin{equation} \label{eqgreen}
\left\{\begin{array}{lc}
\partial_{xx} G_n(x,s)+k_n(x)^2G_n(x,s)=-\delta_s & \text{ in } \R, \\ 
G_n(\cdot,s) \text{ is outgoing.} \end{array}\right. 
\end{equation}
We denote $S_{j_0}$ the region such that $s\in S_{j_0}$, and we assume that $j_0\neq 1$ and $j_0\neq J$ even if it means adding an artificial section before or after the coordinate $s$. Theorem 2 in \cite{olver2} shows that \eqref{eqgreen} has a solution and that there exist $(A_j)_{j=1,...,J},(B_j)_{j=1,...,J}\in \C$ such that $G_n(x,s)$ is close to 
\begin{equation}
G_n^{\text{app}}(x,s)= \left\{\begin{array}{cl} A_1\phi_1(x)w_1(x) & \text{ if } x\in S_1, \\[3pt]
A_J\phi_J(x)w_J(x) & \text{ if } x\in S_J, \\[3pt]
\phi_j(x)(A_j \Ai(\xi_j(x))+B_j\Bi(\xi_j(x))) & \text{ if } x \in S_j, \, j\neq j_0, \\[3pt]
\phi_{j_0}(x)(A_{j_0}\Ai(\xi_{j_0}(x))+B_{j_0}\Bi(\xi_{j_0}(x))) & \text{ if } x\in S_{j_0}, \, x<s, \\[3pt]
\phi_{j_0}(x)(b_\top\Ai(\xi_{j_0}(x))+B_J\Bi(\xi_{j_0}(x))) & \text{ if } x\in S_{j_0}, \, x>s, \end{array}\right. 
\end{equation}
where 
\begin{equation}
w_1(x)=\left\{\begin{array}{cl} \Ai(\xi_1(x)) & \text{ if } h'> 0 \text{ in } S_1, \\ i\Ai(\xi_1(x))+\Bi(\xi_1(x)) & \text{ if } h'< 0 \text{ in } S_1, \end{array}\right.
\end{equation}
\begin{equation} w_J(x)=\left\{\begin{array}{cl} i\Ai(\xi_J(x))+\Bi(\xi_J(x)) & \text{ if } h'>0 \text{ in } S_J, \\ \Ai(\xi_J(x)) & \text{ if } h'<0 \text{ in } S_J. \end{array}\right. 
\end{equation}
To find the value of the constants $A_j$ and $B_j$, we first use the continuity of $G_n^{\text{app}}$ and $\partial_x G_n^{\text{app}}$ on the shared boundaries of each section, which gives $2J-2$ linear equations. Moreover, using the continuity of $G_n^{\text{app}}$ at $x=s$, we find that 
\begin{equation}
A_{j_0}\Ai(\xi_{j_0}(s))+B_{j_0}\Bi(\xi_{j_0}(s))=b_\top\Ai(\xi_{j_0}(s))+B_JB_{j_0}\Bi(\xi_{j_0}(s)).
\end{equation}
Using the jump formula for distributions, we also have 
\begin{equation}
b_\top\Ai'(\xi_{j_0}(s))+B_J\Bi'(\xi_{j_0}(s))-A_{j_0}\Ai'(\xi_{j_0}(s))-B_{j_0}\Bi'(\xi_{j_0}(s))=-\phi_{j_0}(s). 
\end{equation}
Altogether, we obtain a linear system of $2J$ equations for the constants $A_j$, $B_j$, $1\leq j \leq J$. We study its invertibility in the particular case $J=2$ in the next section. 

\subsection{Exemple of dilations or compressions in waveguides}

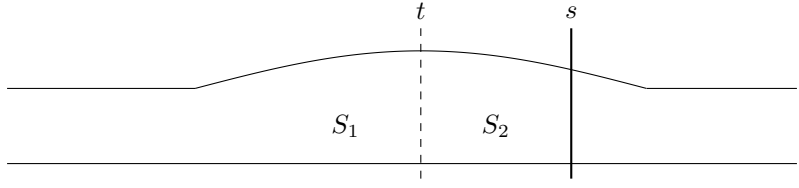
\begin{figure}[h]
\begin{center}
\begin{tikzpicture} 
\draw (-1.5,0) -- (9,0);
\draw (-1.5,1) -- (1,1); 
\draw (7,1) -- (9,1);  
\draw [domain=1:7, samples=200] plot (\x,{0.5*sin(pi*(\x-1)/6 r)+1});  
\draw [dashed] (4, -0.2)--(4,1.8) node[above]{$t$}; 
\draw [thick] (6,-0.2)--(6,1.8) node[above]{$s$}; 
\draw (3,0.5) node{$S_1$}; 
\draw (5,0.5) node{$S_2$}; 
\end{tikzpicture}
\end{center}
\caption{\label{dilatation} Parametrization of a dilated waveguide: $t$ is the unique coordinate in the interior of $\supp h'$ such that $h'(t)=0$, $S_1$ (resp. $S_2$) is the section at the left (resp. right) of $t$, and $s$ is a coordinate satisfying $s>t$.}
\end{figure}

In this section, we apply the method described previously to the simplest case, when the sign of $h'$ changes only once, at $x=t$. If $h$ is increasing then decreasing, we say that the waveguide is dilated. On the other hand, if $h$ is decreasing then increasing, we say that the waveguide is compressed. First, we study the case of dilations. Up to a change of variable $x\mapsto -x$, we can assume that $s>t$, as represented in Figure \ref{dilatation}.

We know from Section 2 that $G_n(x,s)$ is close to 
\begin{equation}\label{greenbis1}
G_n^{\text{app}}(x,s)= \left\{\begin{array}{cl} A_1\phi_1(x)\Ai(\xi_1(x)) & \text{ if } x\in S_1, \\
\phi_{2}(x)(B_1\Ai(\xi_{2}(x))+B_2\Bi(\xi_{2}(x))) & \text{ if } x\in S_{2}, \, x<s, \\
A_2\phi_2(x)\Ai(\xi_2(x)) & \text{ if } x\in S_2 , \,  x>s, \end{array}\right. 
\end{equation}
where $\phi_1$, $\phi_2$ are defined in \eqref{phii}. The constants $A_1,A_2,B_1,B_2$ satisfy the linear system 
\begin{equation} M(s,t)\left(\begin{array}{c}A_1\\B_1\\ B_2\\A_2\end{array}\right)=V(s):=\left(\begin{array}{c} 0\\0\\0\\-\phi_2(s)\end{array}\right),\end{equation}
where $M(s,t)=$
\begin{equation}\label{defm}
\left(\begin{array}{cccc} -(\phi_1\Ai(\xi_1))(t) &(\phi_2\Ai(\xi_2))(t) & (\phi_2\Bi(\xi_2))(t) & 0 \\
-(\phi_1\Ai(\xi_1))'(t) &(\phi_2\Ai(\xi_2))'(t) & (\phi_2\Bi(\xi_2))'(t) & 0 \\
0 & -\Ai(\xi_2(s)) & -\Bi(\xi_2(s)) & \Ai(\xi_2(s)) \\ 
0 & -(\Ai(\xi_2))'(s) & -(\Bi(\xi_2))'(s) & (\Ai(\xi_2))'(s) \end{array}\right).
\end{equation}
We next study when $M$ is invertible. 

\begin{prop}
The determinant $D$ of $M(s,t)$ defined in \eqref{defm} is 
\begin{equation}
D=\frac{1}{\pi}\left(-\phi_1(t)\Ai(\xi_1(t))(\phi_2\Ai(\xi_2))'(t)+(\phi_1\Ai(\xi_1))'(t)\phi_2(t)\Ai(\xi_2(t))\right).\end{equation}
\end{prop}

\begin{proof}
We expand of the determinant along the first column and use the fact that $\Bi\Ai'-\Bi'\Ai=-1/\pi$.
\end{proof}

\begin{rem}
Using the asymptotics of $\Ai$ provided in \cite{perel1, abramowitz1}, the condition $D=0$ asymptotically reduces to 
\begin{equation}
\cos\left(\int_{x^\star_1}^{x^\star_2} k_n(x)\dd x\right)=0.
\end{equation}
Under this condition, it is not possible to find the values of $(A_1,B_1,B_2,A_2)$. This may be explained by the potential presence of trapped modes in the dilated waveguide under this condition. 
\end{rem}

Except for special values of $k$ such that $D=0$, we can find constants $(A_1,B_1,B_2,A_2)$ by computing either symbolically or numerically the solution of $M(s,t)X=V(s)$. 

\vspace{5mm}

We also study the case of compressions. Again, we assume that $s>t$, and we know that $G_n(x,s)$ is close to 
\begin{equation}\label{greenbis2}
G_n^{\text{app}}(x,s)= \left\{\begin{array}{cl} A_1\phi_1(x)(i\Ai(\xi_1(x))+\Bi(\xi_1(x))) & \text{ if } x\in S_1, \\
\phi_{2}(x)(B_1\Ai(\xi_{2}(x))+B_2\Bi(\xi_{2}(x))) & \text{ if } x\in S_{2}, \, x<s, \\
A_2\phi_2(x)(i\Ai(\xi_2(x))+\Bi(\xi_2(x))) & \text{ if } x\in S_2 , \,  x>s. \end{array}\right. 
\end{equation}
The constants $A_1,A_2,B_1,B_2$ satisfy the linear system 
\begin{equation} M(s,t)\left(\begin{array}{c}A_1\\B_1\\B_2\\A_2\end{array}\right)=V(s):=\left(\begin{array}{c} 0\\0\\0\\\phi_2(s)\end{array}\right),\end{equation}
where we define $\mathcal{F}=i\Ai+\Bi$ and $M(s,t)=$
\begin{equation}\label{defm2}
\left(\begin{array}{cccc} -(\phi_1\mathcal{F}( \xi_1))(t) &(\phi_2\Ai(\xi_2))(t) & (\phi_2\Bi(\xi_2))(t) & 0 \\
-(\phi_1\mathcal{F}(\xi_1))'(t) &(\phi_2\Ai(\xi_2))'(t) & (\phi_2\Bi(\xi_2))'(t) & 0 \\
0 & -\Ai(\xi_2(s)) & -\Bi(\xi_2(s)) & \mathcal{F}(\xi_2(s)) \\ 
0 & -(\Ai(\xi_2))'(s) & -(\Bi(\xi_2))'(s) & (\mathcal{F}(\xi_2))'(s) \end{array}\right).
\end{equation}

\begin{prop}
The matrice $M(s,t)$ defined in \eqref{defm2} is invertible. 
\end{prop}

\begin{proof}
Its determinant $D$ is 
\begin{multline}\nonumber
\hspace{1cm}D=\frac{1}{\pi}(-\phi_1(t)(i\Ai+\Bi)(\xi_1(t))(\phi_2(i\Ai(\xi_2)+\Bi(\xi_2))'(t)\\+(\phi_1(i\Ai(\xi_1)+\Bi(\xi_1))'(t)\phi_2(t)(i\Ai+\Bi)(\xi_2(t)).\hspace{3cm}
\end{multline}
The asymptotic expansions provided in \cite{abramowitz1,perel1} show that there exists $\beta_1>0$ such that 
\begin{multline}\nonumber
\phi_1(t)(i\Ai(\xi_1(t))+\Bi(\xi_1(t))\approx\\ \left\{\begin{array}{cl} \frac{1}{2\sqrt{\pi}\sqrt{k_n(t)}}\exp\left({\frac{2}{3}\int_{x^\star_1}^t |k_n|(x)\dd x }\right) & \text{ if } |t-x^\star_1|\gg \eta^{-1/3}, \\ \frac{1}{\sqrt{\beta_1}}
(i\Ai+\Bi)(\beta_1(t-x^\star_1)) & \text{ if } |t-x^\star_1|\ll \eta^{-1/2}, \end{array}\right. 
\end{multline}
and since $h'(t)=0$, 
\begin{multline}\nonumber
(\phi_1(i\Ai(\xi_1)+\Bi(\xi_1))'(t)\approx\\ \left\{\begin{array}{cl} \frac{|k_n(t)|}{2\sqrt{\pi}\sqrt{k_n(t)}}\exp\left({\frac{2}{3}\int_{x^\star_1}^t |k_n|(x)\dd x }\right) & \text{ if } |t-x^\star_1|\gg \eta^{-1/3}, \\ \sqrt{\beta_1}
(i\Ai'+\Bi')(\beta_1(t-x^\star_1)) & \text{ if } |t-x^\star_1|\ll \eta^{-1/2}. \end{array}\right. 
\end{multline}
This case distinction covers every relative positions of $t$, $x^\star_1$ and $x^\star_2$ since $\eta^{-1/2}\gg \eta^{-1/3}$. We can do the same in $S_2$, and find $\beta_2>0$ such that 
\begin{multline}\nonumber
\phi_2(t)(i\Ai(\xi_2(t))+\Bi(\xi_2(t))\approx\\ \left\{\begin{array}{cl} \frac{1}{2\sqrt{\pi}\sqrt{k_n(t)}}\exp\left({\frac{2}{3}\int_t^{x^\star_2} |k_n|(x)\dd x }\right) & \text{ if } |t-x^\star_2|\gg \eta^{-1/3}, \\ \frac{1}{\sqrt{\beta_2}}
(i\Ai+\Bi)(\beta_2(x^\star_2-t)) & \text{ if } |t-x^\star_2|\ll \eta^{-1/2},\end{array}\right. 
\end{multline}
\begin{multline}\nonumber
(\phi_2(i\Ai(\xi_2)+\Bi(\xi_2))'(t)\approx\\ \left\{\begin{array}{cl} -\frac{|k_n|}{2\sqrt{\pi}\sqrt{k_n(t)}}\exp\left({\frac{2}{3}\int_t^{x^\star_2} |k_n|(x)\dd x }\right) & \text{ if } |t-x^\star_2|\gg \eta^{-1/3}, \\ -\sqrt{\beta_2}
(i\Ai'+\Bi')(\beta_2(x^\star_2-t)) & \text{ if } |t-x^\star_2|\ll \eta^{-1/2}.\end{array}\right. 
\end{multline}
If $|t-x^\star_2|\gg \eta^{-1/3}$ and $|t-x^\star_1|\gg \eta^{-1/3}$ then 
\begin{equation}\nonumber
D\approx \exp\left({\frac{2}{3}\int_{x^\star_1}^t |k_n|(x)\dd x }\right)\exp\left({\frac{2}{3}\int_t^{x^\star_2} |k_n|(x)\dd x }\right)\neq 0.
\end{equation}
The remaining  three relative positions of $t$, $x^\star_1$ and $x^\star_2$ can be analyzed in the same way. It follows that $\forall (t,s) \in \R$, $D\neq 0$. 
\end{proof}

Contrarily to the case of dilations, the matrix $M$ is always invertible, and we can compute again the values of $(A_1,B_1,B_2,A_2)$. Figure \ref{figbump} and \ref{fighollow} show examples of such computation.

\section{Numerical illustrations}
In this section, we illustrate our results. We compare the asymptotic expression of $u$ to data generated using the software Matlab to solve numerically the equation \eqref{MAFT} satisfied by the wave field in $\tilde{\Omega}$. In the following, we assume that $h'$ is supported between $x=-7$ and $x=7$. To generate the solution $\wt u$ of \eqref{MAFT} on $\tilde{\Omega}_7$, we use the finite element method and a perfectly matched layer (see \cite{berenger1}) placed on the left side of the waveguide between $x=-15$ and $x=-8$, and on the right side between $x=8$ and $x=15$. The coefficient of absorption for the perfectly matched layer is defined by $\alpha=-k((x-8)\textbf{1}_{x\geq 8}-(x+8)\textbf{1}_{x\leq -8})$ and $k^2$ is replaced in the Helmholtz equation by $k^2+i\alpha$. The structured mesh is built with a stepsize of $10^{-3}$. 

\subsection{Computation of the modal Green function}

\begin{figure}[h!]
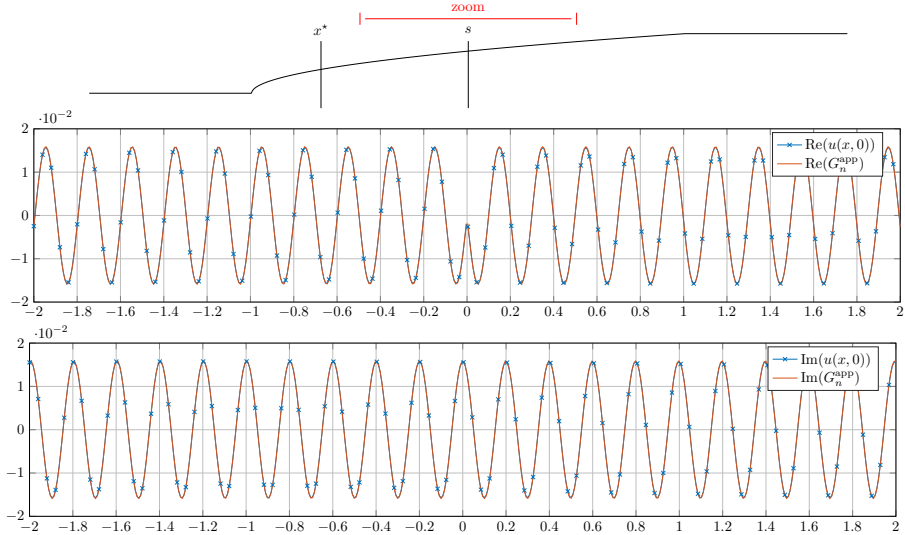

\centering
\hspace{3mm}\scalebox{.55}{\input{h1bis}}
\scalebox{.55}{\input{mode0real}}
\scalebox{.55}{\input{mode0imag}}
\caption{\label{fignleqN} Representation of $\wt u(x,0)$ (in blue) and $G_n^{\text{app}}(x)$ (in red) for a propagative mode $n=0$ and $s=0$ in an expanding waveguide, in the zoomed area $(-2,2)$. Top: representation of $h$, $x^\star$ and $s$. Middle: comparison between real parts of $\wt u(x,0)$ and $G_n^{\text{app}}(x)$. Bottom: comparison between imaginary parts of $\wt u(x,0)$ and $G_n^{\text{app}}(x)$. Here, the relative $\text{L}^2$ error between $\wt u(x,0)$ and $G_n^{\text{app}}(x)$ is $0.93 \%$. }
\end{figure}

To test the validity of the expression \eqref{greenfunction} of the Green function, we consider a profile with $\alpha=0.1$, $\beta=0.04/30$ and 
\begin{equation}\label{hcroiss}
h(x)= \alpha+\beta\left[\left(-1+\frac{\sqrt{x+4}}{{\sqrt{2}}}\right)\textbf{1}_{[-4,4]}(x)+\textbf{1}_{(4,+\infty)}(x)-\textbf{1}_{(-\infty,-4)}(x)\right],
\end{equation} 
with $k=31.5$. There is only one locally resonant mode $N=1$ associated to the resonant point $x^\star\approx -2.72$. We place an internal source $f(x,y)=d(x)\ph_n(y/h(x))$ with $n\in \N$ where $d$ is a Gaussian approximation of $\delta_s$ at $s\in \R$ with
\begin{equation}\label{d}
d(x)=\frac{1}{\sqrt{2\pi}\sigma}\exp\left(-\frac{(x-s)^2}{2\sigma^2}\right),
\end{equation}
and $\sigma=0.005$. We measure the wave field $\wt u$ at $y=0$, and we compare it to the expression \eqref{greenfunction}. In Figure \ref{fignleqN} (resp \ref{figngeqN}), we illustrate the case where $n=0<N$ (resp. $n=2>N$). In Figures \ref{figN1}, \ref{figN2} and \ref{figN3}, we illustrate the case where $n=N$ for different values of $s$. The approximation seems to be accurate, and the small discrepancies observed for instance in the imaginary part of Figure \ref{figN2} are caused by the imprecise approximation of the Dirac function $\delta_s$ by $d$. However, even in this particular case, the relative $\tL^2$ error is still very small. 

\begin{figure}[h!]
\centering
\hspace{1mm}\scalebox{.55}{\input{h1}}
\scalebox{.55}{\input{mode2real}}
\scalebox{.55}{\input{mode2imag}}
\caption{\label{figngeqN} Representation of $\wt u(x,0)$ (in blue) and $G_n^{\text{app}}(x)$ (in red) for an evanescent mode $n=2$ and $s=0$ in an expanding waveguide. Top: representation of $h$, $x^\star$ and $s$. Middle: comparison between real parts of $\wt u(x,0)$ and $G_n^{\text{app}}(x)$. Bottom: comparison between imaginary parts of $\wt u(x,0)$ and $G_n^{\text{app}}(x)$. Here, the relative $\text{L}^2$ error between $\wt u(x,0)$ and $G_n^{\text{app}}(x)$ is $3.68\%$. }
\end{figure}

\begin{figure}[h!]
\centering
\hspace{2mm}\scalebox{.55}{\input{h2}}
\scalebox{.55}{\input{mode11real}}
\scalebox{.55}{\input{mode11imag}}
\caption{\label{figN1} Representation of $\wt u(x,0)$ (in blue) and $G_n^{\text{app}}(x)$ (in red) for a locally resonant mode $n=1$ and $s=-4$ in an expanding waveguide. Top: representation of $h$, $x^\star$ and $s$. Middle: comparison between real parts of $\wt u(x,0)$ and $G_n^{\text{app}}(x)$. Bottom: comparison between imaginary parts of $\wt u(x,0)$ and $G_n^{\text{app}}(x)$. Here, the relative $\text{L}^2$ error between $\wt u(x,0)$ and $G_n^{\text{app}}(x)$ is $5.32\%$. }
\end{figure}

\begin{figure}[h!]
\centering
\hspace{3mm}\scalebox{.55}{\input{h3}}
\scalebox{.55}{\input{mode12real}}
\scalebox{.55}{\input{mode12imag}}
\caption{\label{figN2} Representation of $\wt u(x,0)$ (in blue) and $G_n^{\text{app}}(x)$ (in red) for a locally resonant mode $n=1$ and $s=-1.5$ in an expanding waveguide. Top: representation of $h$, $x^\star$ and $s$. Middle: comparison between real parts of $\wt u(x,0)$ and $G_n^{\text{app}}(x)$. Bottom: comparison between imaginary parts of $\wt u(x,0)$ and $G_n^{\text{app}}(x)$. Here, the relative $\text{L}^2$ error between $\wt u(x,0)$ and $G_n^{\text{app}}(x)$ is $7.24\%$. }
\end{figure}

\begin{figure}[h!]
\centering
\hspace{3mm}\scalebox{.53}{\input{h4}}
\scalebox{.53}{\input{mode13real}}
\scalebox{.53}{\input{mode13imag}}
\caption{\label{figN3} Representation of $\wt u(x,0)$ (in blue) and $G_n^{\text{app}}(x)$ (in red) for a locally resonant mode $n=1$ and $s=5$ in an expanding waveguide. Top: representation of $h$, $x^\star$ and $s$. Middle: comparison between real parts of $\wt u(x,0)$ and $G_n^{\text{app}}(x)$. Bottom: comparison between imaginary parts of $\wt u(x,0)$ and $G_n^{\text{app}}(x)$. Here, the relative $\text{L}^2$ error between $\wt u(x,0)$ and $G_n^{\text{app}}(x)$ is $6.08\%$. }
\end{figure}

We also consider a more general waveguide. As in section 3.2, we choose the simplest case of dilation and compression, and we compare $\wt u(x,0)$ to $G_n^{\text{app}}$ defined in \eqref{greenbis1} and \eqref{greenbis2}. First, we choose to work with a dilated waveguide, described by its width 
\begin{equation}
h(x)=0.1+0.0025\sin\left(\frac{\pi}{10}(x+5)\right)\textbf{1}_{[-5,5]}(x), 
\end{equation}
at frequency $k=31$. The only locally resonant mode is still $N=1$, associated to two resonant points $x^\star_1\approx -3.19$ and $x^\star_2=-x^\star_1$. We choose the same internal source $f$ as before, and we illustrate the case $n=N$ in Figure \ref{figbump}. Figure \ref{fighollow} illustrates the case of a compressed waveguide, with profile
\begin{equation}
h(x)= 0.1-0.0005(x+5)\textbf{1}_{[-5,0]}(x)+\frac{0.0025}{4}(x-4)\textbf{1}_{(0,4]}(x),
\end{equation}
at frequency $k=32.1$ with a resonant mode $N=1$ and resonant points $x^\star_1\approx -0.74$ and $x^\star_2\approx 0.59$. 

\begin{figure}[h!]
\centering
\hspace{3mm}\scalebox{.55}{\input{h5}}
\scalebox{.55}{\input{bumpreal}}
\scalebox{.55}{\input{bumpimag}}
\caption{\label{figbump} Representation of $\wt u(x,0)$ (in blue) and $G_n^{\text{app}}(x)$ (in red) for $n=1$ and $s=0.5$ in a dilated waveguide. Top: representation of $h$, $x^\star_1$, $x^\star_2$ and $s$. Middle: comparison between real parts of $\wt u(x,0)$ and $G_n^{\text{app}}(x)$. Bottom: comparison between imaginary parts of $\wt u(x,0)$ and $G_n^{\text{app}}(x)$ Here, the relative $\text{L}^2$ error between $\wt u(x,0)$ and $G_n^{\text{app}}(x)$ is $10.09\%$. }
\end{figure}

\begin{figure}[h!]
\centering
\hspace{4mm}\scalebox{.53}{\input{h6}}
\scalebox{.53}{\input{hollowreal}}
\scalebox{.53}{\input{hollowimag}}
\caption{\label{fighollow} Representation of $\wt u(x,0)$ (in blue) and $G_n^{\text{app}}(x)$ (in red) for $n=1$ and $s=0.5$ in a compressed waveguide. Top: representation of $h$, $x^\star_1$, $x^\star_2$ and $s$. Middle: comparison between real parts of $\wt u(x,0)$ and $G_n^{\text{app}}(x)$. Bottom: comparison between imaginary parts of $\wt u(x,0)$ and $G_n^{\text{app}}(x)$. Here, the relative $\text{L}^2$ error between $\wt u(x,0)$ and $G_n^{\text{app}}(x)$ is $4.39\%$. }
\end{figure}

\subsection{General source terms}

\begin{figure}[h!]
\centering
\begin{tikzpicture}[scale=1]
\begin{axis}[width=9cm, height=2cm, axis on top, scale only axis, xmin=-7, xmax=7, ymin=0, ymax=0.1013, colorbar,point meta min=1.2606e-05,point meta max=0.0141, title={$|\wt u|$},axis line style={draw=none},tick style={draw=none},scaled y ticks=base 10:2,colorbar style={scaled y ticks=base 10:3}]
\addplot graphics [xmin=-7,xmax=7,ymin=0,ymax=0.1013]{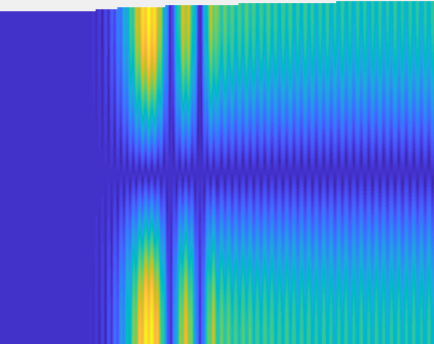};
\end{axis} 
\end{tikzpicture}
\begin{tikzpicture}[scale=1]
\begin{axis}[width=9cm, height=2cm, axis on top, scale only axis, xmin=-7, xmax=7, ymin=0, ymax=0.1013, colorbar,point meta min=2.1962e-07,point meta max=6.3646e-04, title={$|\wt u-\wt u^{\text{app}}|$},axis line style={draw=none},scaled y ticks=base 10:2,tick style={draw=none}]
\addplot graphics [xmin=-7,xmax=7,ymin=0,ymax=0.1013]{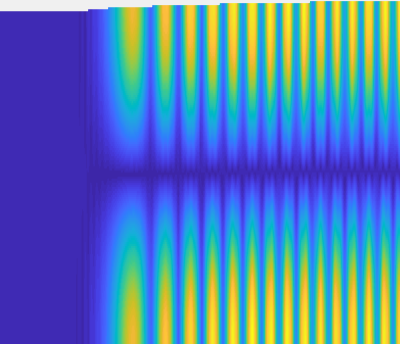};
\end{axis} 
\end{tikzpicture}
\caption{\label{source1} Representation of the wave field $u$ generated by a source $f(x,y)=d(x)y$ where $d$ is defined in \eqref{d} with $s=0$. Up, the absolute value $|\wt u|$, down, the error of approximation $|\wt u-\wt u^{\text{app}}|$. Here, the relative $\text{L}^2$ error between $\wt u$ and $\wt u^{\text{app}}$ is $6.22\%$.}
\end{figure}

\begin{figure}[h!]
\centering
\begin{tikzpicture}[scale=1]
\begin{axis}[width=9cm, height=2cm, axis on top, scale only axis, xmin=-7, xmax=7, ymin=0, ymax=0.1013, colorbar,point meta min=8.6026e-05,point meta max=2.7295, title={$|\wt u|$},axis line style={draw=none},tick style={draw=none},scaled y ticks=base 10:2,colorbar style={scaled y ticks=base 10:1}] 
\addplot graphics [xmin=-7,xmax=7,ymin=0,ymax=0.1013]{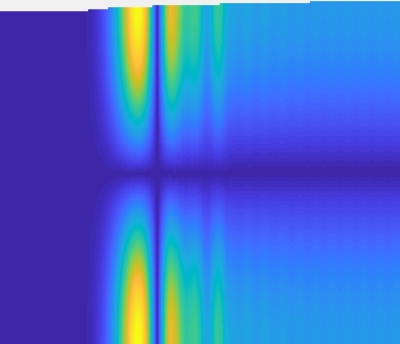};
\end{axis} 
\end{tikzpicture}
\begin{tikzpicture}[scale=1]
\begin{axis}[width=9cm, height=2cm, axis on top, scale only axis, xmin=-7, xmax=7, ymin=0, ymax=0.1013, colorbar,point meta min=1.0695e-04,point meta max=0.2090, title={$|\wt u-\wt u^{\text{app}}|$},axis line style={draw=none},tick style={draw=none},scaled y ticks=base 10:2,colorbar style={scaled y ticks=base 10:2}]
\addplot graphics [xmin=-7,xmax=7,ymin=0,ymax=0.1013]{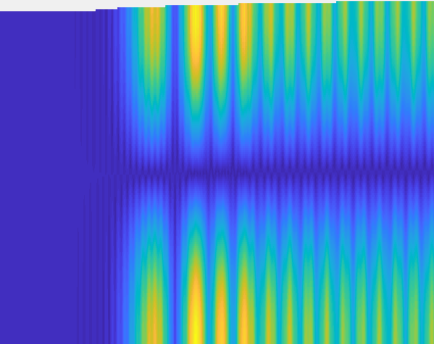}; 
\end{axis} 
\end{tikzpicture}
\caption{\label{source2} Representation of the wave field $\wt u$ generated by a source $b_\bot=\textbf{1}_{-1\leq x\leq 1}$.  Up, the absolute value $|\wt u|$, down, the error of approximation $|\wt u-\wt u^{\text{app}}|$. Here, the relative $\text{L}^2$ error between $\wt u$ and $\wt u^{\text{app}}$ is $7.73\%$.}
\end{figure}

We now validate the approximation provided in \eqref{greentot} for general sources, in the same expanding waveguide defined by its width in \eqref{hcroiss}. We choose two different types of sources :  a vertical internal source (see Figure \ref{source1}), and a boundary source (see Figure \ref{source2}). To compute the approximation in \eqref{greentot}, we choose to reduce the sum to $15$ modes.  Every time, we compute the relative error made between $\wt u$ and its approximation $\wt u^{\text{app}}$ defined by \eqref{greentot}.

\subsection{Dependence of the error of approximation with respect to $\eta$}
Finally, we evaluate in this section the influence of the parameter $\eta$ in the approximation error $\Vert \wt u-\wt u^{\text{app}}\Vert_{\text{H}^1(\widetilde{\Omega}_r)}$. First, we chose large values of $\eta$ to illustrate the control provided in Theorem~\ref{th1}. We work at $k=5.23$, with an increasing waveguide parametrized by 
\begin{equation}
h(x) = 0.5+0.2\eta x\textbf{1}_{-5<x<-5+1/\eta}+0.2\textbf{1}_{x\geq -5+1/\eta},
\end{equation} 
and a boundary source $b_\bot=\textbf{1}_{4<x<5}$. We present in Figure \ref{etagrand} the error of approximation with respect to $\eta$. We can see that the error grows a little bit slower than the $1$ slope expected. We also notice than for small values of $\eta$, the error seems to reach a level where it is almost constant. This is the error due to the finite element method. 

\begin{figure}
\centering
\scalebox{.55}{
%
%
\definecolor{mycolor1}{rgb}{0.00000,0.44700,0.74100}%
\begin{tikzpicture}

\begin{loglogaxis}[%
width=7in,
height=3in,
at={(0in,0in)},
scale only axis,
xmin=0,
xmax=8,
ymin=0,
ymax=300,
grid=minor,
axis background/.style={fill=white},
xlabel={$\eta$},
legend style={legend cell align=left, align=left, draw=white!15!black}
]
\addplot [color=mycolor1]
  table[row sep=crcr]{%
8	43.7205028921113\\
4	30.4464045195003\\
2	29.2790805771293\\
1	18.6144277388636\\
0.333333333333333	10.5365810648273\\
0.222222222222222	7.80534630977656\\
0.166666666666667	5.06083500947094\\
0.125	6.65233026681226\\
0.1	6.0234\\
0.0833333333333333	5.982\\
0.0714285714285714	6.102\\
};
\addlegendentry{$\Vert u-u^{app}\Vert_{\text{H}^1(\widetilde{\Omega}_r)}$}

\addplot [color=black, dashed]
  table[row sep=crcr]{%
8	264.923615669538\\
4	132.461807834769\\
2	66.2309039173846\\
1	33.1154519586923\\
0.333333333333333	11.0384839862308\\
0.222222222222222	7.35898932415385\\
0.166666666666667	5.51924199311538\\
0.125	4.13943149483654\\
0.1	3.31154519586923\\
0.0833333333333333	2.75962099655769\\
0.0714285714285714	2.36538942562088\\
};
\addlegendentry{1 slope}

\end{loglogaxis}

\end{tikzpicture}%
}
\caption{\label{etagrand} Representation of $\Vert \wt u -\wt u^{\text{app}}\Vert_{\text{H}^1(\widetilde{\Omega}_7)}$ for different values of $\eta$ to illustrate the result of Theorem \ref{th1}.}
\end{figure}

\section{Conclusion}

In this paper, we have presented a complete proof of the existence of a unique solution to the Helmholtz equation in slowly variable waveguides. We also provide a suitable approximation of this solution and a control of the error of approximation in $\text{H}^1_\loc(\Omega)$. We validate this approximation numerically, and show that this expression is an excellent way to compute quickly the wave field in a slowly varying waveguide. 

We believe that this work could be extended to elastic waveguides in two dimensions, using the modal decomposition in Lamb modes as in \cite{perel1}. One could also try to generalize the ideas of this article to acoustic waveguides in three dimensions. We think it would be possible to recover some Laplacian eigenvalues of local sections, and from that to recover some information on the waveguide shape. 

Finally, we plan to use this work to develop a new multi-frequency method to recover the width of a waveguide given measurements of the wave field at the surface or on a section of the waveguide. Indeed, for a locally resonant frequency, the wave field in a perturbed waveguide is very different from the one in a regular waveguide, even if the width $h$ is close to a constant function. This should provide a very high sensibility inversion method to reconstruct the width of the waveguide, and will be done in a future work. 

\appendix
\section{Proofs of Proposition \ref{contrs} and \ref{contrb}}

\begin{proof}[Source $f$] 

This proof is an adaptation of the proof presented in Appendix B of \cite{bonnetier1}. Using the results on the modal decomposition presented in Appendix A of \cite{bonnetier1}, we know that the equation \eqref{hdroit} is equivalent to 
\begin{equation}\nonumber
\forall n\in \N \qquad \left\{\begin{array}{cl} u_n''+k_n(x)^2 u_n=-f_n & \text{ in } \R, \\ u_n \text{ is outgoing,} \end{array}\right. 
\end{equation}
where $u_n,f_n\in \text{L}^2_\loc(\R)$ and 
\begin{equation}\nonumber
u(x,y)=\sum_{n\in \N} u_n(x)\ph_n(y), \quad f(x,y)=\sum_{n\in \N} f_n(x)
\end{equation}
Using Theorem \ref{th2}, there exist a unique Green function $G_n(x,s)$ associated to this equation, and
\begin{equation}\nonumber
\forall n \in \N \qquad u_n(x)=\int_\R G_n(x,s) f_n(s) \dd s. 
\end{equation}
We also notice that for every $(x,s)\in \R^2$, $G_n(x,s)=G_n(s,x)$, and using controls from Theorem \ref{th2}, for every $s\in \R$ and $x\in \R$,
\begin{equation}\nonumber
\forall n\in \N \qquad \Vert G_n(\cdot,s)\Vert_{\text{L}^1(-r,r)}, \Vert G_n(x,\cdot)\Vert_{\text{L}^1(-r,r)}\leq \alpha_n^{(1)}.
\end{equation}
Using Young's inequality for integral operators, 
\begin{equation}\nonumber
\forall n\in \N \qquad \Vert u_n\Vert_{\text{L}^2(-r,r)}\leq \alpha_n^{(1)} \Vert f_n \Vert_{\text{L}^2(-r,r)} .
\end{equation}
Using Parseval equality and the results of Theorem \ref{th2},
\begin{multline}\nonumber
\Vert u \Vert^2_{\text{L}^2(\Omega_r)}\leq \alpha^2\sum_{n=0}^{N} \Vert f_n \Vert_{\text{L}^2(-r,r)}^2+\frac{\alpha^2}{\delta^4}\sum_{n>N} \Vert f_n \Vert_{\text{L}^2(-r,r)}^2 \\ \leq \alpha^2\max\left(1, \frac{1}{\delta^4}\right)\Vert f\Vert_{\text{L}^2(\Omega_r)}^2.
\end{multline}
Applying Young's inequality to $u_n'$, we get 
\begin{multline}\nonumber
\Vert \nabla u \Vert^2_{\text{L}^2(\Omega_r)}\leq \alpha^2\sum_{n=0}^{N} \left(1+n^2\pi^2\right)\Vert f_n \Vert_{\text{L}^2(-r,r)}^2\\ +\alpha^2\sum_{n>N} \left(\frac{1}{\min(|k_n|)^2}+\frac{n^2\pi^2}{\min(|k_n|)^4}\right)\Vert f_n \Vert_{\text{L}^2(-r,r)}^2.
\end{multline}
We deduce that 
\begin{equation}\nonumber
\Vert \nabla u \Vert^2_{\text{L}^2(\Omega_r)}\leq \alpha^2\max\left(1+N^2\pi^2,\frac{1}{\delta^2}+\frac{(N+1)^2\pi^2}{\delta^4}\right)\Vert f\Vert_{\text{L}^2(\Omega_r)}^2.
\end{equation}
Finally,
\begin{equation}\nonumber
\Vert u_n''\Vert_{\text{L}^2(-r,r)}\leq |k_n|^2\Vert u_n\Vert_{\text{L}^2(-r,r)}+\Vert f_n \Vert_{\text{L}^2(-r,r)}
\end{equation}
It follows that 
\begin{multline}\nonumber
\Vert \nabla^2 u \Vert_{\text{L}^2(\Omega_r)}^2\leq \sum_{n=0}^{N} \alpha^2\left((k_n^2+1)^2+2n^2\pi^2+n^4\pi^4\right)\Vert f_n\Vert_{\text{L}^2(-r,r)}^2 \\
+\sum_{n>N} \alpha^2\left(\left(\frac{|k_{n}|^2}{\min(|k_{n}|)^2}+1\right)^2+\frac{2n^2\pi^2}{\min(|k_{n}|)^2}+\frac{n^4\pi^4}{\min(|k_{n}|)^4}\right)\Vert f_n\Vert_{\text{L}^2(-r,r)}^2,
\end{multline}
and so 
\begin{multline}\nonumber
\Vert \nabla^2 u \Vert_{\text{L}^2(\Omega_r)}^2\leq  \alpha^2\max\Big[\left((k_N^2+1)^2+2N^2\pi^2+N^4\pi^4\right)\\ +\left(\left(\frac{|k_{N+1}|^2}{\delta^2}+1\right)^2+\frac{2(N+1)^2\pi^2}{\delta^2}+\frac{(N+1)^4\pi^4}{\delta^4}\right)\Big] \Vert f\Vert_{\text{L}^2(\Omega_r)}^2. 
\end{multline}

\end{proof}

\begin{proof}[Source $b$]
Using the same arguments as before, the equation \eqref{hdroitb} is equivalent to
\begin{equation}\nonumber
\forall n\in \N \qquad \left\{\begin{array}{cl} u_n''+k_n^2 u_n=-b_\top\ph_n(1)-b_\bot\ph_n(0) & \text{ in } \R, \\ u_n \text{ is outgoing,} \end{array}\right. 
\end{equation}
and we know that 
\begin{equation}\nonumber
\forall n\in \N \qquad u_n(x)=\int_\R G_n(x,s)(b_\top\ph_n(1)+b_\bot\ph_n(0))\dd s .
\end{equation} 
We notice that $|\ph_n(1)|,|\ph_n(0)|\leq \sqrt{2}$. Using Theorem 2.3.2.9 in \cite{grisvard1}, there exist a constant $d(r)$ and $\mu>0$ such that 
\begin{equation}\nonumber
\Vert u\Vert_{\text{H}^2(\Omega_r)}\leq d(r)\left(\left\Vert -\partial_{xx} u-\frac{1}{h^2}\partial_{yy} u+\mu u\right\Vert_{\text{L}^2(\Omega_r)}+\Vert \bm{b}\Vert_{(\widetilde{\text{H}}^{1/2}(-r,r))^2}\right),
\end{equation}
and it follows that 
\begin{equation}\nonumber
\Vert u\Vert_{\text{H}^2(\Omega_r)}\leq d(r)\left((k^2+\mu)\Vert  u\Vert_{\text{L}^2(\Omega_r)}+\Vert \bm{b}\Vert_{(\widetilde{\text{H}}^{1/2}(-r,r))^2}\right).
\end{equation}
Finally, 
\begin{equation}\nonumber
\Vert u \Vert_{\text{L}^2(\Omega_r)}^2\leq 2\alpha^2\Vert \bm{b}\Vert_{(\widetilde{\text{H}}^{1/2}(-r,r))^2}^2\left(N+1+\sum_{n>N}\frac{1}{\min(|k_n|)^4}\right),
\end{equation} 
which concludes the proof. 
\end{proof}

\bibliographystyle{elsarticle-num}
\bibliography{biblio}

\end{document}